\documentclass[article]{siamart190516}

\usepackage{lipsum}
\usepackage{amsmath}
\usepackage{nccmath}
\usepackage{bm}
\usepackage{amsfonts}
\usepackage{graphicx}
\usepackage{epstopdf}
\usepackage{algorithmic}

\ifpdf
  \DeclareGraphicsExtensions{.eps,.pdf,.png,.jpg}
\else
  \DeclareGraphicsExtensions{.eps}
\fi

\newsiamremark{remark}{Remark}
\newsiamremark{hypothesis}{Hypothesis}
\crefname{hypothesis}{Hypothesis}{Hypotheses}
\newsiamthm{claim}{Claim}

\headers{Decision cycles in a competition}{T.A. McLennan-Smith et al.}

\title{Mathematical modelling of decision cycles in a competitive environment with neutral parties\thanks{\funding{This work was funded by the Defence Science and Technology Group through a Defence Science Partnership agreement}}}

\author{Timothy A. McLennan-Smith \thanks{Defence Science and Technology Group, Canberra, ACT 2600, Australia, (\email{timothy.mclennan-smith@defence.gov.au},\email{alex.kalloniatis@defence.gov.au}).}
\and Alexander C. Kalloniatis \footnotemark[2]
\and Simon Watt \thanks{School of Science, UNSW, Canberra, ACT 2610 , Australia, ( \email{simon.watt@adfa.edu.au}, \email{z.jovanoski@adfa.edu.au}, \email{h.sidhu@adfa.edu.au}, \email{i.towers@adfa.edu.au}).} 
\and Zlatko Jovanoski \footnotemark[3]
\and Harvinder S. Sidhu \footnotemark[3]
\and Dale O. Roberts \thanks{Australian National University, Canberra, ACT 2601, Australia, (\email{dale.roberts@anu.edu.au}).}
\and Isaac N. Towers \footnotemark[3]
} 

\usepackage{amsopn}

\makeatletter
\newcommand*{\addFileDependency}[1]{ 
  \typeout{(#1)}
  \@addtofilelist{#1}
  \IfFileExists{#1}{}{\typeout{No file #1.}}
}
\makeatother

\begin{document}

\maketitle

\begin{abstract}
Decision making is a human process that is a fundamental part of competition. As a realisation of decision making, Command and Control, or C2, has been studied in the literature for adversarial populations, yet these models do not explicitly capture competition. Our work here seeks to enrich such competition models by extending ecologically inspired population models to include decision dynamics through coupling with a nonlinear system of oscillators (the Kuramoto equation) to model the action-perception cycle in decision making.  Through asymmetric competition models of increasing complexity, we highlight the importance of competitive agility and decision making processes necessary for a population to be successful. 
\end{abstract}




\section{Introduction}

Complex systems, by their nature, involve multiple elements interacting with each other, typically across multiple purposes and objectives; the expression `multi-layered networks' is often used to reflect such complexity in modern theories~\cite{Boccaletti2014TheSA}. In that context, we present and study various dynamical models in which competition and cooperation play out, regulated by layered decision-making processes. As an example, consider the model
\begin{align}
    \dot{P}_1 & = r_1 P_1\left(1 - P_1 \right ) - \beta_2 P_1 P_2 \, (\sin(\overline{\theta}_{2}-\overline{\theta}_{1}) +2)/2, \label{eq:modelsimple1}  \\
    \dot{P}_2 & =  r_2  P_2\left(1 - P_2 \right ) - \beta_1 P_2 P_1 \, (\sin(\overline{\theta}_{1}-\overline{\theta}_{2}) +2)/2, \label{eq:modelsimple2}  \\
    \dot{\theta}_i &= \omega_i +  (1 - P_2\mathbf{1}_{i\in V_1}-P_1\mathbf{1}_{i\in V_2})  \, \textstyle\sum_{j\in V} [W_{ij} \sin (\theta_j - \theta_i + \Phi_{ij})], \label{eq:modelsimple3} 
\end{align}
for $i \in V$ where $P_1(t)$ and $P_2(t)$ represent the magnitude of resources of each population. This is a simplified two-population version, of more generalised models that we consider later in this paper, where we have assumed that both populations are synchronised in their decision cycles. Notice that this model combines a predator-prey type Lotka-Volterra model with logistic recruitment through the first term in \eqref{eq:modelsimple1} and \eqref{eq:modelsimple2}, combined with \eqref{eq:modelsimple3} which models a `decision cycle' through the coupling of a Kuramoto-Sakaguchi type model on networks defined by the graphs $G_1 = (V_1,E_1)$ and $G_2 = (V_2,E_2)$ and the combined graph $G = (V,E)$ where $|V| = N$; see~\cite{dorogovtsev2008critical,arenas2008synchronization,dorfler2014synchronization}. The $N \times N$ weighted adjacency matrix $\mathbf{W}$ encodes the presence of edges $E$ between nodes and the magnitude of coupling between them. The phases $(\theta_i)_{i=1}^N \in [0,2\pi)$ are regarded as a `decision states' of agents within a cyclic decision-state cognitive model, in essence, a continuous version of Boyd's Observe-Orient-Decide-Act (OODA) loop~\cite{boyd1987organic,kalloniatis2020EJOR}. In this respect, native frequencies $\omega_i$ represent intrinsic decision making speeds of agents, networks their patterns of interaction, coupling strengths the intensity of that interaction, and `frustration' in the decision-making dynamics through the matrix $\mathbf{\Phi} = (\Phi_{ij})$ which quantifies how far ahead or behind agent $i \in V$ seeks to be in relation to the decision state of agent $j \in V$. Here $r_1$ and $r_2$ represent the intrinsic recruitment/growth rates in the logistic part of the model. The collective decision cycles of the agents  on networks $G_1$ and $G_2$ influence the competitive dynamics through the sine difference of the centroids $\overline{\theta}_{1} := \frac{1}{N_1} \sum_{i\in V_1} \theta_{i}$ and $\overline{\theta}_{2} := \frac{1}{N_2} \sum_{j\in V_2} \theta_{j}$ where $\beta_2$ is the rate at which $P_2$ reduces $P_1$ in competition and $\beta_1$ is the rate at which $P_1$ reduces $P_2$. Thus being ahead in its collective decision state (up to $\pi$) allows a population to more effective in the competition.

We present the model in \eqref{eq:modelsimple1}-\eqref{eq:modelsimple3} and its generalisations framed within a military context of competing and adversarial populations. That said, our approach is also applicable to economic, financial, administrative, and social fields where distributed agents have varying degrees of agency in outcomes of the system behaviour. Moreover, beyond human applications, such dynamical systems can be considered in predator-prey ecologies where animal cognition constitutes the decision-making process~\cite{d2011cognitive}. The aim of our paper is to show how such models may be formulated mathematically and the range of analytical and computational tools that may be brought to bear upon them. Our investigations will take the reader through versions of the model with increasing complexity, where the focus of our analysis will primarily be placed upon the effect of introducing a decision-making layer into the competitive predator-prey style dynamics. We identify distinct regions in select parameter subspaces, where decision making processes produce thresholds for the success for one side in the competition and consider a design of experiments approach to find the key parameters that play the biggest role in controlling the dynamics of the model. This allows us to demonstrate how the choice of the feedback from decision synchronisation onto the competition can affect relative parameter significance. 



\section{Background and motivation}\label{sec:background}

Dynamical system approaches to represent cooperating and competing actor-populations have a long history in the literature~\cite{smith1986asymptotic,cushing1989competition}. More recently, modern network theoretic tools have become prominent in such systems \cite{agbanusi2018emergence}. For decision making processes, we have representations such as the Learning Intelligent Distribution Agent (LIDA) scheme for human cognition \cite{baars2007architectural} or an animal cognitive extension \cite{d2011cognitive}. Within military and business strategy there is the OODA loop model \cite{boyd1987organic}. All of these share the common hypothesis that cognitive processing is via repeated iteration of a fundamental action-perception cycle \cite{Neisser1967,schoner1998action,freeman2002limbic} and it now well understood that these dynamics often appear in Command and Control (C2), a term from military, civilian emergency-response, policing and plant management contexts \cite{athans1987command,holling1996command}. C2 as a distributed decision making system is a fundamentally human process, the domain of organisational theorists, psychologists and business analysts. However efforts to apply mathematics to this problem space have been long-standing \cite{dockery1985mathematics}. Recently, an approach has been advanced employing the dynamical system of the Kuramoto model \cite{kalloniatis2020EJOR} as a representation of interacting agents connected on networks, attempting to synchronise individual action-perception cycles. The Kuramoto model remains an active area of research in the context of synchronisation on networks, \cite{kuehn2019power,skardal2019synchronization} and in describing the structure of ecosystems \cite{vandermeer2021new}.

Earlier applications of the Kuramoto model treated the synchronisation of decision cycles in isolation from an explicit representation of the external purpose of decision-making \cite{kalloniatis2016fixed,holder2017gaussian,kalloniatis2019two}. However, competition is ultimately a contest of rival resource strengths. Within the typical military domain of C2, the Lanchester equations provide a compact mathematical description of competition between adversarial military forces \cite{lanchester1916aircraft,taylor1983lanchester}. Recognising that the Lanchester equations can be considered an instantiation of the Lotka-Volterra model for two species in a predator-prey relationship (the simplest model of an ecological system) provides a straightforward way to generalise to multiple actors in a competitive environment. This is realised in the generalisation beyond `trophic' effects. As increasingly recognised, some species in an ecological network play no trophic role at all - that is they are neither preyed upon, nor prey others. Rather, they reside in niches amongst other foraging species, facilitating the growth of one or another element \cite{kefi2012more,fontaine2011ecological,berlow2004interaction,navarrete1990resource,melian2009diversity,mougi2012diversity,gross2008positive}. Recent work \cite{mclennansmith2021} has generalised the two-party Lanchester model to accommodate such a diverse range of actors. This separate line of research has focused on modelling the interactions of adversaries within a host population environment, also seen in recent conflicts \cite{zuparic2022modelling}. From a military perspective, C2 is a key enabler of kinetic effects as can be observed in contemporary conflict. This pivotal role of decision making through C2 serves as the main motivation for the development of our hybrid model and we seek to address its sparsity in previous Lanchester-type modelling literature. 

A further motivation for this work is to complement the types of models typically used in the contemporary defence context spawned by the temptation of unlimited computational power, often using detailed agent-based distillation or alternate physics motivated representations. Examples of these include the StateSim model of Silverman and collaborators \cite{silverman2021statesim}, at one extreme, and Ising model
inspired representations of societal dynamics
of Axelrod and Bennett \cite{Axelrod1993},
its adaptation to C2
by Song et al. \cite{Song2013}, or the approach of Carley and collaborators \cite{Carley1999}. These
are but a few examples as the literature of these
approaches is large. The key challenge of these
models lies in the significant number of
parameters required, or the complete divorce
from otherwise trusted approaches to conflict
dynamics as the Lanchester system. The modelling
approach espoused in our work, and predecessor
papers, seeks to enhance the traditional Lanchester approach to include decision-making
within a variety of populations, within
a more data-parsimonious framework. 

Validation of all such models, including ours,
is a challenge; a richer suite of models
permits greater cross-validation 
\cite{Bharathy2010} of concepts
and hypotheses across a spectrum of modelling
approaches. Within this, the ability to bring
closed-form analytical solutions alongside
high performance computational approaches
is critical. In this paper then we bring together the decision making representation of Kuramoto dynamics with this multi-population model of trophic and non-trophic effects, building on \cite{ahern2021}. In this sense we fuse into one, several disparate threads of complex networks, competition and ecological modelling within a mathematical representation of decision making to achieve the desired hybrid dynamical system.


\section{Linearisation of the competition model} \label{section:simplemodel}

In order to study the model presented in \eqref{eq:modelsimple1}-\eqref{eq:modelsimple3}, we first consider the following block forms of the weighted adjacency and frustration matrices
\begin{align}
\mathbf{W} = \left( \begin{array}{cc}
     \sigma_1  \mathbf{A}^{(11)} & \xi_{12}  \mathbf{A}^{(12)} \\
   \xi_{21}   \mathbf{A}^{(21)} & \sigma_2   \mathbf{A}^{(22)} 
\end{array}\right), \;    \mathbf{\Phi} = \left( \begin{array}{ccc}
        0 & \phi  \textbf{J}^{(12)}  \\
        \psi  \textbf{J}^{(21)} & 0 \\
    \end{array}\right), \label{eq:matrices_2pop}  
\end{align}
where $\mathbf{A}^{(ii)}$ the adjacency matrix defined by the edges $E_i$, $\mathbf{A}^{(ij)}$ is the adjacency matrix defining by edges between  $V_i$ and $V_j$ in the combined graph $G$, and $\mathbf{J}^{(ij)}$ is matrix of ones equal in size to $\mathbf{A}^{(ij)}$. Using this form of the frustration matrix $\mathbf{\Phi}$, inner-population dynamics will seek to synchronisation their decision cycles whilst inter-population dynamics will attempt to be ahead of their competitor by amounts $\phi$ and $\psi$. From the Kuramoto dynamics, we define the quantities
\begin{align}
\overline{\omega}_i := \frac{1}{N_i} \sum_{i \in V_i} \omega_i , \quad \mu := \overline{\omega}_1 - \overline{\omega}_2, \quad d_k^{(ij)} := \sum_l A^{(ij)}_{kl}, \label{eq:networkproperties_2pop} 
\end{align}
where $\overline{\omega}_i$ is the average native frequency value of nodes $V_i$, $\mu$ is the difference in average native frequencies between $V_1$ and $V_2$, and $d_k^{(ij)}$ is the degree of $V_k$ in graph $G_i$ to graph $G_j$. We assume that both populations are synchronised in their decision cycles. This assumption allows for application of two previous results taken from earlier works  \cite{kalloniatis2016fixed} and \cite{kalloniatis2020manoeuvre}, with derivations provided in the Supplementary material. The first of these results provides a first-order reduction of the phase dynamics through the following proposition:

\begin{proposition} \label{theorem:twocluster}
Given two internally synchronised populations with phase dynamics described by \eqref{eq:modelsimple3} with the centroid difference 
$\Delta = \overline{\theta}_{1}-\overline{\theta}_{2}$,
then a first-order approximation of the phase dynamics in \eqref{eq:modelsimple3} is given by 
\begin{align}
    \dot{\Delta} = \mu + S \cos \Delta - C \sin \Delta,
\end{align}
where
\begin{align}
 C & \equiv (d_T^{(12)} \xi_{12} H_1/ N_1 ) \cos \phi  + (d_T^{(21)} \xi_{21} H_2/N_2)\cos \psi , \\
  S & \equiv (d_T^{(12)} \xi_{12} H_1/N_1 ) \sin \phi - (d_T^{(21)} \xi_{21} H_2/N_2)\sin \psi , 
\end{align}
with $H_1 = 1-P_2$, $H_2 = 1-P_1$, and $d_T^{(\cdot)} = \sum_{i \in V_1} d_{i}^{(\cdot)}$. 
\end{proposition}
\begin{proof}
The proof follows from a previous result in \cite{kalloniatis2020manoeuvre}. 
\end{proof}

Thus we can apply Proposition~\ref{theorem:twocluster} to approximate the set of ODEs in \eqref{eq:modelsimple3} to the following single ODE for the centroid difference:
\begin{align}
 \dot{\Delta}  = \mu + S \cos \Delta - C \sin \Delta, \label{eq:simplecentroidODE}
\end{align}
where $\Delta := \overline{\theta}_{1}-\overline{\theta}_{2}$, and we can directly substitute $\Delta$ into \eqref{eq:modelsimple1} \& \eqref{eq:modelsimple2}.

We now note that for some fixed $\Delta^* \in \mathbb{R}$, the fixed points of the system, denoted by $\text{FP}_i =(P^*_1,P^*_2 ,\Delta^*)$, may be solved giving:
\begin{align}
    \text{FP}_1 & = (1,0,\Delta^*), \quad \text{FP}_2 = (0,1,\Delta^*),\quad \text{FP}_3 = (0,0,\Delta^*),\\
    \text{FP}_4 & =  \left(2\frac{r_2 \left( 2r_1 + \beta_2 \sin(\Delta^*) -2 \beta_2 \right) }{4r_1r_2 +\beta_1 \beta_2 \left(\sin^2(\Delta^*) -4\right)},\, 2\frac{r_1\left( 2r_2 - \beta_1 \sin(\Delta^*) -2 \beta_1 \right) }{4r_1r_2 +\beta_1 \beta_2 \left(\sin^2(\Delta^*) -4\right)},\Delta^*\right). \label{eq:simpleE4}
\end{align}
In order to handle the $\Delta^*$ term, we observe the second result from \cite{kalloniatis2016fixed} and \cite{kalloniatis2020manoeuvre} which provides a solution to the centroid dynamics $\dot{\Delta}$ through the following proposition:

\begin{proposition} \label{prop:centroiddynamics}

For some fixed $P^*_1,P^*_2 \in [0,1]$, the solution to the centroid dynamics for two internally synchronised populations is given by
\begin{align} 
\Delta(t) =  2 \tan^{-1} \left(  \left( C-\sqrt{\mathcal{K} }\tanh ( \sqrt{ \mathcal{K}} (t+\textrm{const.} )/2  ) \right) \middle/ \left(\mu-S\right) \right),
\label{eq:Deltasol}
\end{align}
where 
\begin{eqnarray}\label{eq:special-K}
\mathcal{K} \equiv C^2 + S^2 -\mu^2. 
\end{eqnarray}
Furthermore, the existence of a fixed point requires that 
\begin{align}
    C^2 + S^2 \geq \mu^2.
\end{align}
\end{proposition}
\begin{proof}
The proof follows from a previous result in \cite{kalloniatis2020manoeuvre}. 
\end{proof}

\begin{figure}
    \centering
    \includegraphics[width=\textwidth]{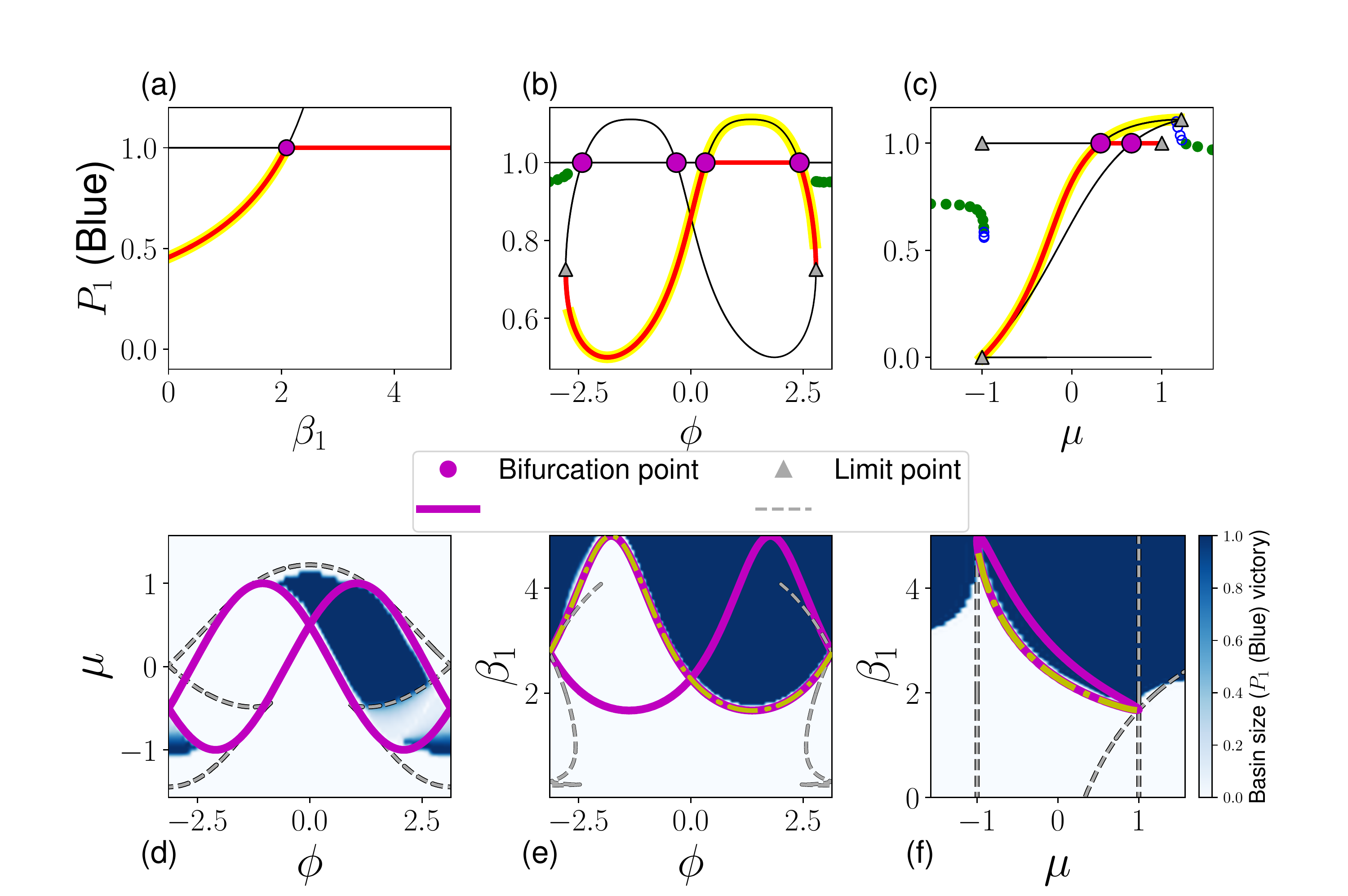}
    \caption{Single parameter (a)-(c) and two parameter (d)-(f) bifurcation diagrams for the two population model with linearisation from \eqref{eq:simplecentroidODE}. In the single parameter plots, red lines denote stable fixed points, black lines denote unstable fixed points, open blue circles denote unstable limit cycles, and solid green circles denote stable limit cycles. Yellow lines denote the fixed points obtained analytically from \eqref{eq:system_simple}. For the two parameter (d)-(f) bifurcation diagrams, we include the basin of attraction size for Blue success. The yellow dotted-dashed lines denotes the bound from \eqref{eq:betaCondition}.
    (a): $\phi = 0.2$, $\mu = 0.2$, (b): $\beta_1 = 2$, $\mu = 0.2$, (c): $\beta_1 = 2$, $\phi = 0.2$, (d): $\beta_1 = 2$, (e): $\mu = 0.2$, and (f) $\phi = 0.2$. }
    \label{fig:simplemodel}
\end{figure}

From Proposition~\ref{prop:centroiddynamics}, we take the following result for the fixed point of the centroid difference, given some $P_1^*,P_2^* \in [0,1]$, through the limit
\begin{align} \label{eq:delta0}
    \Delta^* = \lim_{t\rightarrow \infty} \Delta(t) =  2 \tan^{-1} \left(  (C-\sqrt{\mathcal{K} })/(\mu-S ) \right).
\end{align}
For the first three fixed points where there is no $\Delta^*$ dependency in $P_1^*$ or $P_2^*$, we directly substitute the fixed points $P_1^*$ and $P_2^*$ into \eqref{eq:delta0} to obtain
\begin{align*}
\begin{split}
    \text{FP}_1:& \quad \Delta^* = 2 \tan ^{-1}\left(  (\gamma_1 \cos (\phi )-\sqrt{\gamma_1^2 -\mu ^2})\middle/(\mu -\gamma _1 \sin (\phi ))\right), \\ 
    \text{FP}_2:& \quad  \Delta^* = 2 \tan ^{-1}\left((\gamma_2 \cos (\psi )-\sqrt{\gamma_2^2 -\mu ^2})\middle/(\mu +\gamma_2 \sin (\psi ))\right),\\ 
    \text{FP}_3:& \quad \Delta^* = 2 \tan ^{-1}\left(\frac{\gamma_1 \cos (\phi )+\gamma_2 \cos (\psi ) - \sqrt{\gamma_1^2 + \gamma_2^2 +2\gamma_1\gamma_2 \cos (\phi +\psi ) -\mu ^2}}{\mu -\gamma _1 \sin (\phi )+\gamma_2 \sin (\psi )} \right),
\end{split}\label{eq:fpsimpledelta}
\end{align*}
where $\gamma_1 = {d_T^{(12)} \xi_{12} }/{N_1}$ and $\gamma_2 = {d_T^{(21)} \xi_{21} }/{N_2}$. For the fixed point $\text{FP}_4$ in \ref{eq:simpleE4} which contains a $\Delta^*$ dependency for both $P_1^*$ and $P_2^*$, the fixed point may be explicitly obtained by solving the following system of equations through fixed point iteration
\begin{align}
\begin{split}
    P_1^* &  =2\frac{r_2 \left( 2r_1 + \beta_2 \sin(\Delta^*) -2 \beta_2 \right) }{4r_1r_2 +\beta_1 \beta_2 \left(\sin^2(\Delta^*) -4\right)}, \quad
    P_2^* = 2\frac{r_1\left( 2r_2 - \beta_1 \sin(\Delta^*) -2 \beta_1 \right) }{4r_1r_2 +\beta_1 \beta_2 \left(\sin^2(\Delta^*) -4\right)}, \\
\Delta^* &=  2 \tan^{-1} \left( \frac{C(P_1^*,P_2^*)-\sqrt{\mathcal{K}(P_1^*,P_2^*) }}{\mu-S(P_1^*,P_2^*)}\right).
\end{split} \label{eq:system_simple}
\end{align}
Solutions to this system are compared to the fixed points obtained by numerical continuation in Figure~\ref{fig:simplemodel}. Before discussing the plot in detail, we test the stability properties of the fixed points by examining the Jacobian matrix
{\small
\[\begin{pmatrix}
 (2 P_1-1) r_1+\frac{1}{2} \beta_2  P_2 (\sin \Delta-2) & \frac{1}{2} P_1 \beta_2  (\sin \Delta-2) & \frac{1}{2} \beta_2 P_1 P_2 \cos \Delta \\
 -\frac{1}{2} \beta_1  P_2 (\sin \Delta+2) & r_2 (1-2 P_2)-\frac{1}{2} \beta_1  P_1 (\sin \Delta+2) & -\frac{1}{2} \beta_1  P_1 P_2 \cos \Delta \\
\gamma _2 \sin (\Delta+\psi) & \gamma _1 \sin ( \Delta - \phi ) & -S \sin \Delta - C \cos \Delta  \\
\end{pmatrix}\!.
\]
}%
For the fixed points $\text{FP}_1,\text{FP}_2$, and $\text{FP}_3$, the Jacobian yields the following eigenvalues whose sign will be analysed for stability/instability
\begin{align*}
    \text{FP}_1: & \quad  \lambda_{11}  = - r_1,  \; \lambda_{12}  = \tfrac{1}{2}(-2 \beta_1 - \beta_1 \sin \Delta^* + 2r_2), \;  \lambda_{13}  = -\gamma_1\cos(\phi - \Delta^*). \\
    \text{FP}_2: & \quad   \lambda_{21}  = \tfrac{1}{2}(-2 \beta_2 + \beta_2 \sin \Delta^* + 2r_1), \; \lambda_{22}  = -r_2, \; \lambda_{23}  =  -\gamma_2\cos(\psi + \Delta^*) . \\
    \text{FP}_3: &  \quad   \lambda_{31}  =  r_1, \; \lambda_{32} = r_2, \; \lambda_{33}  = -\gamma_1\cos(\phi - \Delta^*) -\gamma_2\cos(\psi + \Delta^*).
\end{align*}

Observing the eigenvalues $\lambda_{12}$ and $\lambda_{22}$ from the fixed points $\text{FP}_1$ and $\text{FP}_2$, where $\text{FP}_1$ and $\text{FP}_2$ corresponds to $P_1$ and $P_2$ victories respectively, we obtain the stability conditions:
\begin{align}
   \lambda_{12}: \quad \beta_1 \, > \,& \frac{r_2}{1+\frac{1}{2} \sin(\Delta^*)}; \quad
   -\frac{\pi}{2} +\pi k \, \leq \, \phi-\Delta^* \, \leq \, \frac{\pi}{2} + \pi k, \quad k \in \mathbb{Z}, \label{eq:betaCondition}\\
   \lambda_{21}: \quad  \beta_2 \, > \,& \frac{r_1}{1-\frac{1}{2} \sin(\Delta^*)}; \quad
   -\frac{\pi}{2} +\pi k \, \leq \, \psi+\Delta^* \, \leq \, \frac{\pi}{2} + \pi k, \quad k \in \mathbb{Z}. \label{eq:etaCondition}
\end{align}
The trivial fixed point $\text{FP}_3$ will always be unstable since $r_1,r_2 >0$ by definition.

We now consider a specific instance of the two population competition model in \eqref{eq:modelsimple1}-\eqref{eq:modelsimple3}  where the phase dynamics from \eqref{eq:modelsimple3} are reduced to \eqref{eq:simplecentroidODE}. In this case study, we fix $\gamma_1 = 1$, $\gamma_2= 1$, $\psi=0$, $\beta_2=2$, $r_1 =3$ and $r_2 =2.5$. For the intrinsic frequencies of the populations $\omega_i$, we choose randomly from the uniform distribution, that is $\omega_{i} \sim U[0,1]$ for $i \in V_1,V_2$. We select this limited range to represent the scenario that the populations are internally closely aligned in their intrinsic decision making as opposed to $[0,2\pi]$ which represents full diversity of decision making. In this scenario, where in typical Lanchester fashion we name $P_1$ as Blue and $P_2$ as Red, both sides interact with each other with the same coupling strength, however Blue is given a slightly greater logistic capability. Critically, Red in this setting is not seeking to be ahead in decision cycles through $\psi =0$. A final differentiation is taken by setting the graph $G_1$ to a $k$-ary tree and $G_2$ to be an  Erd\"os-R\'enyi graph. We will continue with this setup later in Section~\ref{sec:networks}.   
We provide a more detailed discussion of the numerical simulation in the Supplement which includes our routine for the calculation of $\overline{\theta}$ as directly taking the arithmetic mean for the centroid can yield incorrect results if the winding number $q$ is non-zero. As part of our numerical results, we present a bifurcation analysis of the system obtained using the software XPPAUT \cite{ermentrout2002simulating}, where we focus upon $\beta_1,\mu$ and $\phi$ as the bifurcation parameters interpreting these as the variables Blue controls (its ability to reduce Red, its average decision making speed, and how much it aspires to be ahead of Red decision making) for achieving success.

As a metric for parameter regimes that lead to Blue success, we consider the proportional size of the basin of attraction of the fixed point $E_1$ that corresponds to Blue success or survival, where the optimal value of this proportion is 1 meaning that Blue will always survive (and Red extinguished) regardless of initial conditions. In the calculation of the basin of attraction, we limit the phase space of initial conditions between 0 and 1 (or the carrying capacity in the dimensional case) for each population. 

Our numerical investigation, together with the key analytical insights from this section, are summarised on Figure~\ref{fig:simplemodel}. Blue is able to survive with more than half its initial resources (a stable fixed point) for reduction rates less than that of Red $\beta_2=2$, as a consequence of its faster decision speeds and decision advantage with $\mu=\phi=0.2$ -- seen in panel (a).  When Blue's reduction rate is equal to Red's, its success is guaranteed even when it has decision disadvantage (negative $\phi$) as long as it has faster decision making speeds, evident in the range of the stable fixed point in panel (b). In panel (c) we see that, with equal reduction rate to Red, and with Blue exercising decision advantage $\phi>0$, Blue shows a stable fixed point with corresponding surviving populations even when it has a range of slower decision speeds with $\mu<0$. We also note the important role of the sign of $\mathcal{K}$ from \eqref{eq:delta0}. With the $\mu$ bifurcation plot in (c) as an example, we observe that $\mathcal{K} < 0 $ for $|\mu| > 1$, thus the existence of our fixed points is constrained to the space $|\mu| \leq 1$. For $|\mu| > 1$ we observe the appearance of limit cycles. In summary, plots such as these show how the different mechanisms available to Blue for enhancing its position combine to guarantee its success even when disadvantaged in other aspects. 

\section{The generalised competition model}\label{sec:fullmodel}
Expanding out from the simple two population competition model in  \eqref{eq:modelsimple1}-\eqref{eq:modelsimple3}, we now outline the general form of the model treated in this paper. 

\subsection{Modelling framework}
Let $P_{i} \in [0,K_{i}]$, denote resources of the $i^{th}$ population group with carrying capacity $K_i$ where $i \in [1,N]$. Let $G_{j} = (V_{J},E_{j})$ denote the graph nodes of population decision agents of size $N_j$ and $G = (V,E)$ denote the combined population graph.  The dynamics of the competition between the collective populations is given by the system of ordinary differential equations written as
\begin{align}
    \dot{P}_i & = L_i  p(\mathcal{O}_{S_{i}}) - \sum_{j \in [1,N]}F_{ji}  P_j  I(\mathcal{O}_{T_{j}},\bar{\theta}_i,\bar{\theta}_j) - X_{i}, \quad  i \in [1,N], \label{eq:combatmodel_L} \\
    \dot{\theta}_k &= \omega_k + H_k  \sum_{l\in V} W_{kl}  \sin (\theta_l - \theta_k + \Phi_{kl}) , \quad k \in V.
\label{eq:combatmodel_K}
\end{align}
In this hybird system, \eqref{eq:combatmodel_L} represents the predator-prey (Lanchester) dynamics between the competing population groups and \eqref{eq:combatmodel_K} describes the decision dynamics (Kuramoto-Sakaguchi) across the population groups. Under this general formulation, the population resource dynamics take inspiration from Lotka-Volterra type models where we incorporate three commonly occurring effects of growth, competition and decay through the functions $L_i$, $F_{ij}$, and $X_i$ respectively. The first of these functions $L_i := L_i(P_1,P_2,\dots)$ encodes population logistics/growth which are specified for each population. The competition functions $F_{ji}:= F_{ji}(P_i)$ denote rates of reduction of population $P_i$ by population $P_j$; these may also be described as the reductive power or strength of one group over the other and arise as the result of direct competition. The functions $X_i := X_i(P_1,P_2,\dots)$ denotes competition fatigue or withdrawal. In the case where $\forall i \in [1,N]: F_{ki} = 0 $, we regard population $P_k$ as not directly competing in the given scenario whilst a population $P_l$ is considered to actively compete if $ \exists i\in [1,N]:  F_{li} > 0$. Although a population $P_k$ might not explicitly compete, it may still influence the Lanchester dynamics of other populations $P_i$ by indirect competition (or support) through $L_i$ and $X_i$ and may implicitly modulate other $F_{ij}$ terms through decision making feedback from the function $I$. Internal decision sychronisation is described through the function $p(\mathcal{O})$ where  $\mathcal{O} \in [0,1]$ is the Kuramoto order parameter \cite{kuramoto2003chemical}. 
Synchronisation of decision making within a population is thus regarded as an enabler of logistics and competition as it directly modulates the functions $L_i$ and $F_{ij}$ through the order parameter. In applying the order parameter we consider the order parameters $\mathcal{O}_{S_{i}}$ and $\mathcal{O}_{T_{i}}$ which represent `strategic' $(S_i)$ and `tactical' $(T_i)$ sub-populations within the population $V_i$ given by
\begin{align}
    \mathcal{O}_{i} := \frac{1}{N_{i}} \left| \sum_{j \in V_{i} } e^{j \theta_{j}} \right|, \quad  \mathcal{O}_{S_{i}} := \frac{1}{N_{S_i}} \left| \sum_{j \in S_i } e^{j \theta_{j}} \right|, \quad      \mathcal{O}_{T_{i}} := \frac{1}{N_{T_i}} \left| \sum_{j \in T_i } e^{j \theta_{i}} \right|. 
\end{align}
While strategic is suggestive of a senior leadership, executive, or species apex role, the term tactical may be interchanged with `worker', `client facing', `front-line' or `foragers' depending on the application context. Fundamentally, some form of social stratification is assumed in such structures which distinguishes how different members interact with the competitor population. We however note this stratification is scenario dependent and that these sub-groups may be simplified by instead considering the overall population synchronisation $\mathcal{O}_{i}$.  

The feedback function $I := I(\mathcal{O}_{T_{j}},\bar{\theta}_i,\bar{\theta}_j)$ is regarded as the initiative function which enhances the respective reduction rate $F_{ji}$ based upon the collective decision states of the populations $i$ and $j$ through the phase centroids $\bar{\theta}_i$ and $\bar{\theta}_j$ respectively as well as the sychronisation of competing agents with population $j$ through $\mathcal{O}_{T_{j}}$. Generally this function will enhance competition function $F_{ij}$ if a population $P_i$ is ahead of $P_j$ in its collective decision state. Finally feedback from the Lanchester dynamics onto the Kuramoto coupling is described through the function $H_k: = H_k(P_1,P_2,\dots)$ which scales the ability of the $k^{th}$ agent within $V$ to synchronisation with other agents as a function of population resources from the Lanchester dynamics.   

\subsection{The two population model with synchronisation feedback}

We now revisit the previous two population competition \eqref{eq:modelsimple1}-\eqref{eq:modelsimple3} with the inclusion of synchronisation feedback through the following system

\begin{align}
    \dot{P}_1 & = r_1 P_1\left(1 - P_1 \right )  \mathcal{O}_{S_{1}}  - \beta_2 P_1 P_2  \mathcal{O}_{T_{2}} (\sin(\overline{\theta}_{2}-\overline{\theta}_{1}) +2)/2, \label{eq:modelsimpleex1}  \\
    \dot{P}_2 & =  r_2  P_2\left(1 - P_2 \right )  \mathcal{O}_{S_{2}} - \beta_1 P_2 P_1  \mathcal{O}_{T_{1}} (\sin(\overline{\theta}_{1}-\overline{\theta}_{2}) +2)/2, \label{eq:modelsimpleex2}  \\
    \dot{\theta}_i &= \omega_i +  (1 - P_2\mathbf{1}_{i\in V_1}-P_1\mathbf{1}_{i\in V_2})   \sum_{j\in V} W_{ij} \sin (\theta_j - \theta_i + \Phi_{ij}), \quad i \in V, \label{eq:modelsimpleex3} 
\end{align}
We note that this system is an instance of the general model \eqref{eq:combatmodel_L}-\eqref{eq:combatmodel_K} with the following functional forms in the Lanchester dynamics
\begin{align}
    L_i &= r_i P_i(1-P_i), \quad F_{ji} = (1-\delta_{ij})\beta_i P_i, \quad  X_i = 0, \quad i,j \in [1,2].
\end{align}
In coupling these models together we employ feedbacks between decision making and the competition, described through the linear function $p(\mathcal{O}_{(\cdot)}) = \mathcal{O_{(\cdot)}}$ and
\begin{align} \label{eq:feedbackforms}
    I(\mathcal{O}_{T_i},\bar{\theta}_j,\bar{\theta}_k) = \mathcal{O}_{T_i}  \frac{\sin(\bar{\theta}_k - \bar{\theta}_j) +2}{2}, \quad 
    H_i(P_1,P_2) = \left\{ \begin{array}{cc}
         1 - P_2 & \text{if } i \in V_1 \\[2mm]
         1 - P_1 & \text{if } i \in V_2\\[2mm] 
    \end{array}
    \right. .
\end{align}
This choice of $I$ is based on the difference in the population phase centroids $\bar{\theta}_1,\bar{\theta}_2$ which represents how far ahead a competing population is in its average decision cycle relative to its adversary and further modulates the ability of a population to engage in the competition based on the synchronisation of its tactical agents through $\mathcal{O}_{T_i}$. Feedback of the competition on the decision coupling through the function $H$ is based upon the assumption that a larger adversarial population inhibits the effectiveness of a member of the competing population to interact with others. We thus consider $H$ as a measure of stress that is placed upon collective decision capabilities as a result of adversarial competition, where at full stress saturation members will resort to decision making at their own intrinsic rates. 

\begin{figure}
    \centering
    \includegraphics[width=\textwidth]{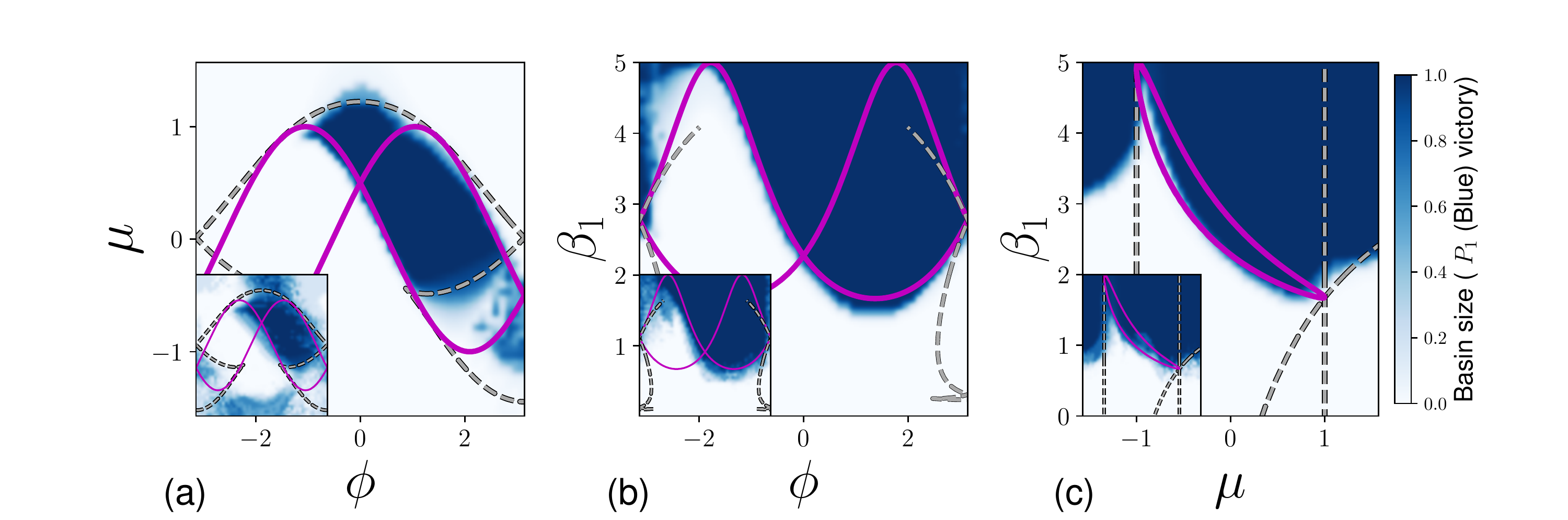}
    \caption{Basin of attraction size for Blue success for the model in \eqref{eq:modelsimpleex1}-\eqref{eq:modelsimpleex3} with $\sigma_1=4, \sigma_2 = 2$ (high synchronisation) and (a): $\beta_1 = 2$, (b): $\mu = 0.2$, and (c) $\phi = 0.2$. Bifurcation lines from Figure~\ref{fig:simplemodel} are superimposed for comparison. Inset plots provide an alternative parameterisation of the system with weaker coupling values of $\sigma_1=1, \sigma_2 = 0.5$ (low synchronisation).}
    \label{fig:simplemodel2}
\end{figure}

For this version of the two population competition model with synchronisation feedback, we set $\sigma_1=4$ and $\sigma_2 = 2$ and set $\xi_{12} = N_1/d_T^{(12)}$  and $\xi_{21} = N_2/d_T^{(21)}$ to match the choices of $\gamma_1 = 1$ and $\gamma_2= 1$ from the linearised case study of \eqref{eq:modelsimple1}-\eqref{eq:modelsimple3} in Section~\ref{section:simplemodel}. These represent relatively high values of coupling for the types of networks used, promoting, though not necessarily guaranteeing, synchronisation given the forms of non-linear feedback in the model. In addition, we consider an additional parameterisation where $\sigma_1 = 1$ and  $\sigma_2 = 0.5$ to highlight the divergence in behaviour when synchronisation is less likely. We plot these two cases of high and low synchronisation across panels (a)-(c) in Figure~\ref{fig:simplemodel2} with the latter as an inset plot. For the high synchronisation case, we see close agreement on the form of the basin of attraction of $\text{FP}_1$ (Blue victory) with that from the linearisation results shown in Figure~\ref{fig:simplemodel}. To aid in this comparison we have superimposed the bifurcation lines from Figure~\ref{fig:simplemodel}.  We observe indeed that the synchronised approximation in \eqref{eq:simplecentroidODE} produces similar behaviour to the synchronisation version in \eqref{eq:modelsimpleex1}-\eqref{eq:modelsimpleex3} with agreement on stable fixed points across broad parameter ranges when the model is parameterised to promote sychronisation amongst the populations. In contrast, the parameterisation with weaker coupling values, which results in low sychronisation, begins to diverge from the structures from Figure~\ref{fig:simplemodel} yet retains qualitative similarities. 

\section{An ecologically inspired case of the competition model} \label{section:ecomodel}

We now consider a three population case of the competition based on ecological modelling through the presence
of a third, non-trophically interacting, population \cite{kefi2012more} and the functional response \cite{holling1959some}. The version of the extended model that we consider generalises our previous work \cite{mclennansmith2021} to incorporate decision making across the three {\it networked} population dynamics. As in the previous scenario, there exists two populations $P_1$ (Blue) and $P_2$ (Red) in competition. In the conflict interpretation of this model, Green, or $P_3$, may be a group providing aid or refuge to one or the other side in a purely altruistic mode, but relies on some form of protection or infrastructure support by the other side. As Green is a non-competitor within our modelling framework as expressed via \eqref{eq:combatmodel_L}-\eqref{eq:combatmodel_K}, this requires that $\forall i \in [1,3]: F_{3i} = 0$. However Green will influence the Lanchester dynamics of the other two competitors in \eqref{eq:combatmodel_L} through non-trophic interactions.

\subsection{Network structure selection}\label{sec:networks}

\begin{figure}
    \centering
    \includegraphics[width=1.\textwidth]{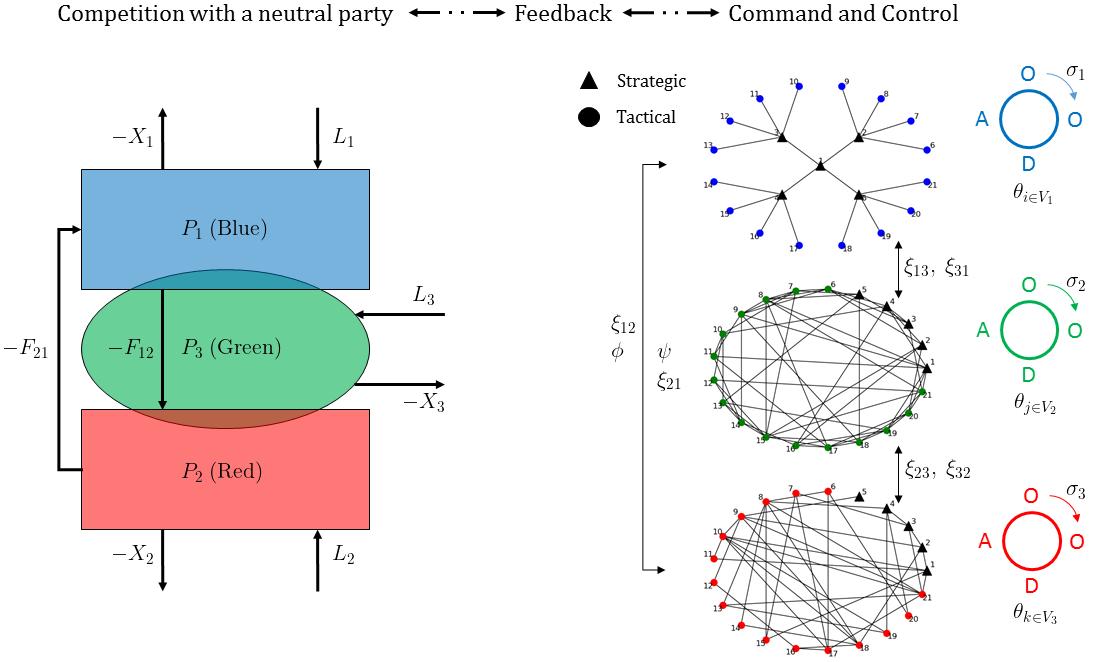}
    \caption{Interaction diagram of the decision making competition model with a neutral party and non-trophic interactions. Left: Solid black lines represent the various interactions of recruitment, competition reduction, and fatigue/withdrawal present in the model given by \eqref{eq:combatmodel_L}. Right: Sample networks and decision dynamics described through Kuramoto phase oscillators.}
    \label{fig:interactions}
\end{figure}

In this paper, we consider a particular use case for the networks for Blue, Green and Red that represent stylised forms similar in size but distinct in organisational structure for the three populations in consideration. 

Consistent with previous work \cite{kalloniatis2016fixed,holder2017gaussian,kalloniatis2019two,mclennansmith2021} where we contrast strongly different forms of organisation, we treat Blue strictly as hierarchical in its communication/interaction pathways, thus use a $k$-ary tree as our graph structure with four branches and two layers; the graph size of Blue is $N_1 = 21$. For Red, we use a random graph to represent an organisation without hierarchical structure in its communication paths, though still with a distinction between strategic/apex and tactical/foragers. Here the random graph is an Erd\"os-R\'enyi graph  \cite{erdos2011evolution} with size $N_2 = 21$ and a link probability of $0.2$ between nodes. Finally, for Green we use the small-world network \cite{watts1998collective} typical of social structures \cite{hong2002synchronization} of size $N_3 = 21$ with a rewiring probability of $0.4$, where each node is initially connected to 6 nearest neighbours in a ring topology. On Figure~\ref{fig:interactions} we present the graphs and their connecting edges for this use case. For the Blue nodes labelled 1-5, there are connecting edges to the Green nodes labelled 1-5. For the Blue nodes labelled 6-21, there are connecting edges to the Red nodes labelled 6-21. Finally, for the Red nodes labelled 6-21, there are connecting edges to the Green nodes labelled 6-21. Note that the Red nodes labelled 1-5 have no inter-network connections. For these three populations, we consider the following extension of the weighted adjacency and frustration matrices from Eq.\eqref{eq:matrices_2pop} as follows
\begin{align}\label{eq:matrices}
W = \left( \begin{array}{ccc}
     \sigma_1  \mathcal{A}^{(11)} & \xi_{12}  \mathcal{A}^{(12)} & \xi_{13}  \mathcal{A}^{(13)} \\
   \xi_{21}   \mathcal{A}^{(21)} & \sigma_2   \mathcal{A}^{(22)} & \xi_{23}  \mathcal{A}^{(23)}\\
    \xi_{31}  \mathcal{A}^{(31)} & \xi_{32}  \mathcal{A}^{(32)} & \sigma_3  \mathcal{A}^{(33)}
\end{array}\right), \;    \Phi = \left( \begin{array}{ccc}
        0 & \phi  \textbf{J}^{(12)}  & 0 \\
        \psi  \textbf{J}^{(21)}  & 0 & 0 \\
        0& 0 & 0
    \end{array}\right),
\end{align}
where $\mathbf{A}^{(ij)}$ and $\mathbf{J}^{(ij)}$ are the adjacency matrices and matrices of ones as previous discussed. For the intrinsic frequencies of the populations $\omega_i$, as before we select $\omega_{i} \sim U[0,1]$ for $i \in V_1,V_2$ where $U$ is the uniform distribution. For Green, we set $\omega_i = 0.5$ for $i \in V_3$ as Green is to be viewed as a tight-knit community \cite{newman2000models}. 

\subsection{Model formulation}

For the three population scenario with two competitors (Blue and Red) and one non-competitor (Green), we consider the following model:
\begin{align}
\begin{split}
    \dot{P}_1 & = r_1^*  P_1(1-{P_1}/{K_1}) \mathcal{O}_{S_1} - F^*_{21} P_2 \mathcal{O}_{T_2} ({\sin(\bar{\theta}_{2} - \bar{\theta}_{1}) +2})/{2} - x^*_{1} P_1, \\
    \dot{P}_2 & =  r_2^* P_2(1-P_2/K_2) \mathcal{O}_{S_2} - F^*_{12} P_1 \mathcal{O}_{T_1} (\sin(\bar{\theta}_{1} - \bar{\theta}_{2}) +2)/{2} - x^*_{2} P_2 , \\
    \dot{P}_3 & =  r_3^*   P_3(1-{P_3}/{K_3}) \mathcal{O}_{S_3} -x^*_3 P_3,  \\
    \dot{\theta}_i &= \omega_i + H_i \sum_{j \in V}W_{ij} \sin (\theta_j - \theta_i + \Phi_{ij}) , \quad i \in V.
    \end{split} \label{eq:fullcombatmodel}
\end{align}
The first two equations here represent the predator-prey (Lanchester) dynamics between the two competing groups, Blue and Red; the third equation brings the neutral population `Green' into play; and the final equation gives the overall decision dynamics across the three groups. This is an instance of the general model \eqref{eq:combatmodel_L}-\eqref{eq:combatmodel_K} with the following functional forms in the Lanchester dynamics

\begin{align}
    L_i &= r^*_i P_i(1-P_i), \quad X_i = x_i^*P_i, \quad i,j \in [1,2,3], \\
    F_{ji} &= \left\{ \begin{array}{ll}                                         F^*_{ji} & \text{ if } j=1,i=2 \text{ or } j=2,i=1, \\
                                        0 & \text{ otherwise}.
    \end{array} \right. 
\end{align}
As before, we employ feedback between decision making and the competition, described through the linear function $p(x) = x$ and
\begin{align} \label{eq:feedbackforms_eco}
    I(\mathcal{O}_{T_i},\bar{\theta}_j,\bar{\theta}_k) = \mathcal{O}_{T_i} \frac{\sin(\bar{\theta}_k - \bar{\theta}_j) +2}{2}, \quad 
    H_i = \left\{ \begin{array}{cc}
         1 - P_2 & \text{if } i \in V_1, \\[2mm]
         1 - P_1 & \text{if } i \in V_2,\\[2mm] 
         1  & \text{if } i \in V_3.\\[2mm] 
    \end{array}
    \right. 
\end{align}
This is the hybrid dynamical system we study in the rest of this paper and now proceed to describe all the terms present in this model. In this system the non-trophic inspired interactions are given by 
\begin{align}
    r^*_1 &= r_1 \frac{ \alpha  P_2}{1+ \alpha  P_2}, \quad r^*_2= r_2,  \quad r^*_3= \frac{r_3 +r^{(\max)}_{3}  P_1}{1+P_1},  \label{eq:nontrophic1}  \\ F^*_{12} &= \frac{ \beta_1^* P_2 }{1+  \tau \beta_1^* P_2  }, \quad
\beta_1^*  = \frac{\beta_1+  \beta_1^{(\min)} P_3  }{1+ P_3}, \quad F^*_{21} = \beta_2 P_1, \label{eq:nontrophic2} \\ x^*_1 &= x_1, \quad x^*_2 = 0, \quad
    x^*_{3} = x_3 -\frac{ \left(x_3 - x^{(\min)}_{3}\right) P_1 }{1+P_1 } + \frac{ \left(x^{(\max)}_{3} - x_3\right) P_2 }{1+P_2 }. \label{eq:nontrophic3}
\end{align}
Here recruitment occurs through a logistic growth term with a carrying capacity $K$ and an intrinsic recruitment rate $r$  for each population. In this formulation, we introduce a parameter $\alpha$ now to reflect how quickly Blue is able to
shift its resources in response to the presence of Red.  In the Lanchester setting of the model \cite{syms2015dynamic}
this reflects the real-world limitations in an organisation assembling and maintaining resources from a point external to the environment where the adversary may be found. In this scenario Blue will only recruit if it detects the presence of Red. Effectively we use a sigmoidal response to switch on Blue recruitment where the sharpness of this response is described by the parameter $\alpha$, which we denote as the strategic agility of Blue. Thus, the larger the value of $\alpha$, the faster Blue is able to react to the presence or emergence of Red. The rate of Green logistic recruitment can be enhanced to $r^{(\max)}_{3}$ by the presence of Blue, as Blue supports a secure environment where Green can perform its refuge provisioning role. Competition between Blue and Red occurs in an asymmetric manner. Red upon Blue competition is described in a Lotka-Volterra type interaction of $-\beta_2 P_1 P_2$ where $\beta_2$ is the rate that Red reduces Blue in competition.

\begin{figure}
    \centering
    \includegraphics[width=\textwidth]{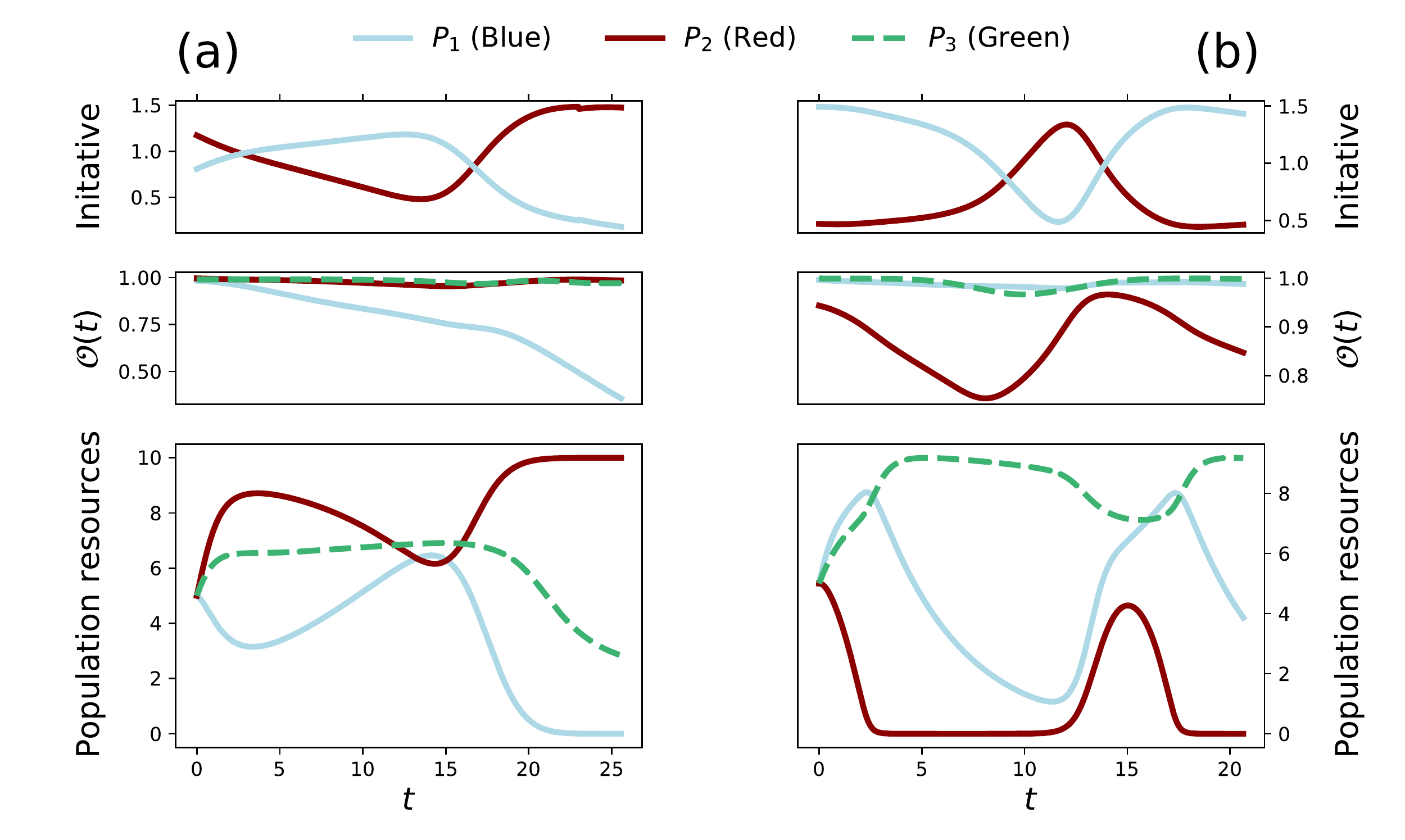}
    \caption{Sample dynamics for two cases of the full model in \eqref{eq:fullcombatmodel} with (a) $\beta_1=1.5$, $\phi=-\pi/2$, and (b) $\beta_1=5$, $\phi=\pi/2$. For each case, time series plots display population resources, the order parameter $\mathcal{O}(t)$ for the whole populations, and the competition initiative for the Blue and Red populations.}
    \label{fig:TSplot}
\end{figure}

Note how we have switched the use of the symbol $\beta_1$ in Eq.\eqref{eq:nontrophic2} compared with its use as
Blue's rate of reducing Red in the simpler model Eq.\eqref{eq:modelsimpleex1}. The parameter $\beta_1$ will
still describe the competitive capability of Blue but now, in the spirit of using ecological analogues, 
we use it in the construct of the Holling functional response \cite{holling1959some} which we denote by $F^*_{12}$. In ecological modelling, this describes the response intake rate of a consumer as a function of food density. In our context, 
reflecting some of the realities of contemporary conflict environments as explained in \cite{mclennansmith2021},
we employ a switching between an `indirect' to `direct' reduction rate based upon the size of Red present. 
The rate of this switching is determined by the expression $\beta_1^*$ which we term the tactical agility of Blue. The intrinsic tactical agility of Blue is given by the parameter $\beta_1$ and in the absence of Green we have $\beta_1^*=\beta_1$. At saturation in the Holling factor, Red's reduction follows a `direct' form and becomes $-P_2/\tau$ where $\tau$ is the searching and engagement time of Blue. In other words, for increasingly large $\beta_1^*$ the functional response is only as good as the time it takes for Blue to find and compete with Red. For small values of Red, reduction approximately follows the `indirect' form of $-\beta_1^* P_1 P_2$ as in the simplified model in \eqref{eq:modelsimpleex2}.

The non-trophic interaction through Green in the functional response occurs through refuge provisioning for Red. It thus modifies the tactical agility of Blue to tend towards $\beta_1^{(\min)}$ which is the restricted tactical agility of Blue in the presence of Green. As a neutral altruistic party, Green does not suffer reduction. For the final fatigue/withdrawal interactions, Blue will withdraw from the competition at a constant rate $-x_1$ which can describes a fatigue and/or the general withdrawal, when it views its competition as a success. We thus regard a Blue success as having had occurred when both Blue and Red resources are 0. Green suffers fatigue $-x_3$ due to being in a competitive environment, in which fatigue can be mitigated due to the supportive efforts of Blue to $x^{(\min)}_{3}$, and exacerbated due to the antagonistic efforts of Red to $x^{(\max)}_{3}$.  We note that in this formulation we set to unity all sigmoidal sharpness parameters involving and affecting Green interactions in the $r^*_3$, $x^*_3$, and $\beta_1^*$ terms.

For our general case study of this model, we use: $\nu = -0.25$, 
$\psi=0$, $\beta_2 =0.2$, $r_1 =3$, $r_2 =2.5$, $r_3=1$, $\alpha=2$,  $\tau = 1$, $x_1=0.25$, $x_3=0.25$, $\beta_1^{(\min)} =0.1$, $x_3^{(\min)}=0.125$, $x_3^{(\max)}=0.5$, $r_3^{(\max)} = 1.5$, and $K_{1}=K_{2}=K_{3}=10$. 
These choices provide paradigmatic examples of the types of dynamics that may occur. In Figure~\ref{fig:TSplot} we show two cases that may be realised where Red succeeds in panel (a) and where Blue succeeds in panel (b). The two cases here differ solely in the choice of the decision parameter $\phi$, where (a) $\phi=-\pi/4$ and (b) $\phi = \pi/4$, thus the cases in (a) and (b) correspond to Blue permitting a decision disadvantage and seeking advantage respectively. In the top panels of these plots we show the time series of the initiative functions $I$ from \eqref{eq:feedbackforms} for both Blue and Red; in the middle panels we give the overall population synchronisation order parameter $\mathcal{O}_i(t)$ which combines both tactical and strategic sub-populations; and on the bottom panels we present the population dynamics for the three populations. We observe here how a combination of higher levels of decision cycle synchronisation and decision advantage determines the outcome. Recalling that initiative acts as a modulation of the reduction rates in the competitions, a loss of Blue initiative drives a resurgence of Red as observed between times $t\in [10,15]$ in a panel (b). Previous work \cite{mclennansmith2021} presented the phenomenon of Red resurgence as a consequence of Blue withdrawal which also persists in this model. However, in the present model the additional attenuation of the reduction rate through decision disadvantage enhances the capability of Red to surge. Ultimately in panel (b), Blue is able to regain decision superiority and through it defeat Red.

\subsection{The centroid approximation}

The results in Propositions~\ref{theorem:twocluster} and~\ref{prop:centroiddynamics} are however only applicable for the two population version of our model. Assuming that all three populations in \eqref{eq:fullcombatmodel} are close to internal synchronisation ($\mathcal{O}_1=\mathcal{O}_2=\mathcal{O}_3=1$), we present the following we consider the following proposition
\begin{proposition} \label{theorem:threecluster} 
Given three internally synchronised populations with phase dynamics described by \eqref{eq:fullcombatmodel} and the centroid differences
\begin{align}
    \Delta_1 = \overline{\theta}_{1} - \overline{\theta}_{2}, \quad \Delta_2 = \overline{\theta}_{1}- \overline{\theta}_{3}, \quad \Delta_3 = \overline{\theta}_{2} - \overline{\theta}_{3},
\end{align}
then the first-order approximation for synchronised phase dynamics in \eqref{eq:fullcombatmodel} are given by
\begin{align}
    \dot{\Delta}_1 & = \mu - H_1 \left( (\xi_{12} d_T^{(12)}/N_1) \sin\left( \Delta_1 -\phi \right) - (\xi_{13} d_T^{(13)}/N_1) \sin\left( \Delta_2 \right) \right)  \\
    & -H_1  \left( (\xi_{21} d_T^{(21)}/N_2) \sin\left( \Delta_1 +\psi \right) + (\xi_{23}  d_T^{(23)}/N_2) \sin\left( \Delta_2 - \Delta_1\right) \right), \nonumber \\
    \dot{\Delta}_2 & = \nu -H_2 \left(  (\xi_{12} d_T^{(12)}/N_1) \sin\left( \Delta_1 -\phi \right) - \xi_{13} d_T^{(13)}/N_1) \sin\left( \Delta_2 \right) \right) \\
    &  - (\xi_{31} d_T^{(31)}/N_3) \sin\left( \Delta_2 \right) - (\xi_{32} d_T^{(32)}/N_3) \sin\left( \Delta_2 - \Delta_1\right). \nonumber
\end{align}
where $\Delta_3 = \Delta_2 - \Delta_1$, $H_1 = (1-P_2)$, and $H_2 = (1-P_1)$.

\end{proposition}

\begin{proof}
The proof follows from a previous result in \cite{kalloniatis2016fixed}. 
\end{proof}

We can apply Proposition~\ref{theorem:threecluster} to reduce the set of ODEs for $\theta$ to the following two ODEs for the two centroid differences, $\Delta_1$ and $\Delta_2$, as follows

\begin{equation}
\begin{aligned}
    \dot{P}_1 & = r^*_1  P_1\left(1 - P_1/K_1 \right ) - \beta_2 P_1 P_2  (\sin(-\Delta_1)+2)/2 - x_1 P_1, \\
    \dot{P}_2 & =  r_2 P_2 \left(1 - P_2/K_2 \right ) - F^*_{12} P_1  (\sin(\Delta_1)+2)/2 , \\
    \dot{P}_3 & = r^*_3  P_3 \left(1 - P_3/K_3  \right ) - x^*_3 P_3, \\
     \dot{\Delta}_1 & = \mu - (1-P_2) \left( (\xi_{12} d_T^{(12)}/N_1) \sin\left( \Delta_1 -\phi \right) + (\xi_{13} d_T^{(13)}/N_1) \sin\left( \Delta_2 \right) \right)  \\
    & - (1-P_1) \left( (\xi_{21} d_T^{(21)}/N_2) \sin\left( \Delta_1 +\psi \right) - (\xi_{23}  d_T^{(23)}/N_2) \sin\left( \Delta_2 - \Delta_1\right) \right), \\
    \dot{\Delta}_2 & = \nu - (1-P_2) \left(  (\xi_{12} d_T^{(12)}/N_1) \sin\left( \Delta_1 -\phi \right) + (\xi_{13} d_T^{(13)}/N_1) \sin\left( \Delta_2 \right) \right) \\
    &  - (\xi_{31} d_T^{(31)}/N_3) \sin\left( \Delta_2 \right) - (\xi_{32} d_T^{(32)}/N_3) \sin\left( \Delta_2 - \Delta_1\right),
\end{aligned} \label{eq:fullmodelapprox}
\end{equation}
where $r_1^*$,$F_{12}^*$, $r_3^*$, and $x_3^*$ are defined in \eqref{eq:nontrophic1}-\eqref{eq:nontrophic2}. 

\begin{figure}
    \centering
    \includegraphics[width=1\textwidth]{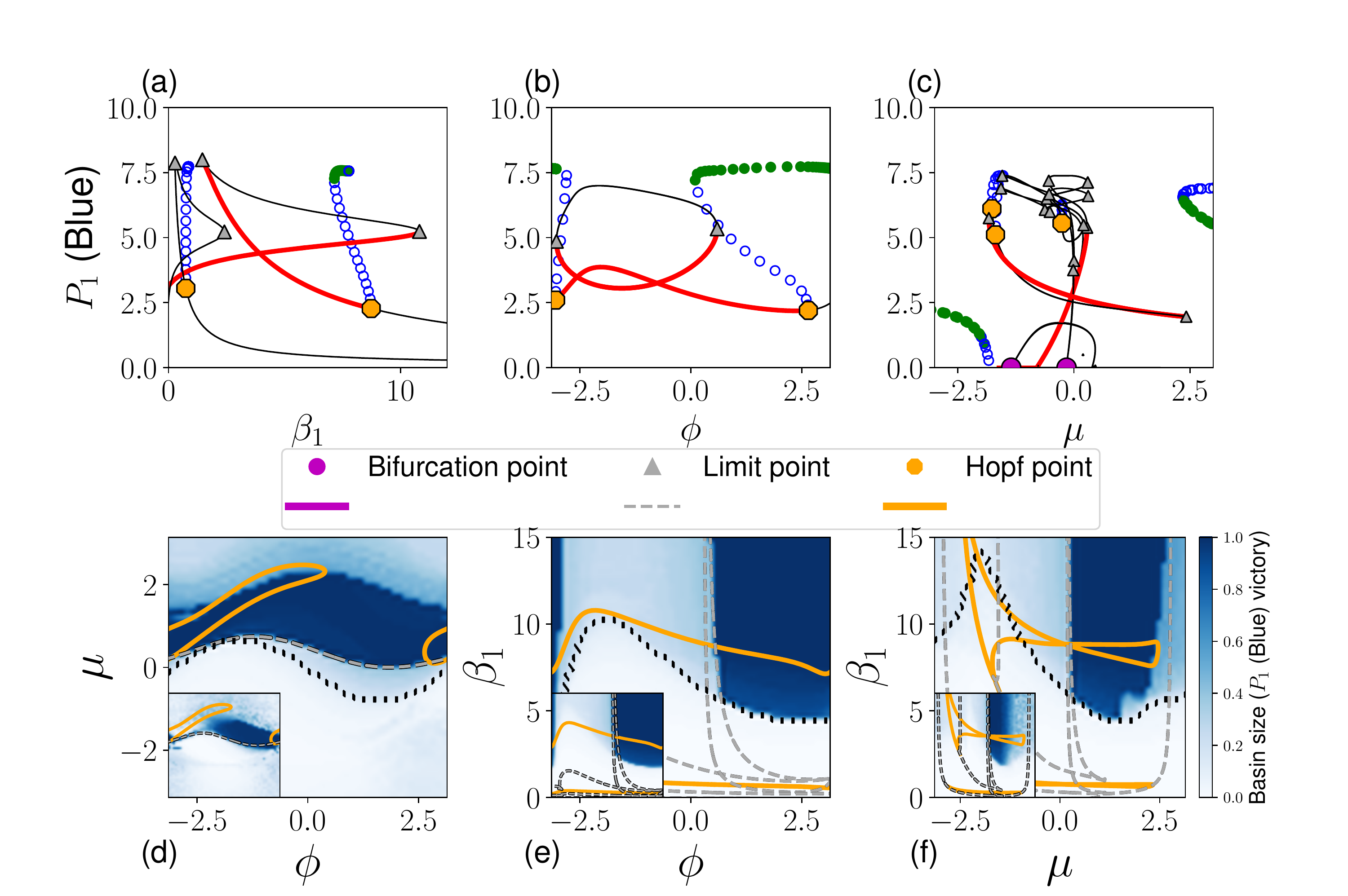}
    \caption{Single parameter (a)-(c) and two parameter (d)-(f) bifurcation diagrams for the linearisation of the competitive Kuramoto model from \eqref{eq:fullmodelapprox} where (a) $\phi = 0.5$, $\mu = 0.25$ (b) $\beta_1 = 7.5$, $\mu = 0.25$, (c) $\beta_1 = 7.5$, $\phi = 0.5$, (d) $\beta_1 = 7.5$, (e) $\mu = 0.25$, (f) $\phi = 0.5$.  Here bifurcation lines and basin values are denoted as in Figure~\ref{fig:simplemodel}. Dotted black lines represent the contour of the basin values above 0.15. Inset plots show the basins of attraction for the networked model in \eqref{eq:fullcombatmodel} with $\sigma_1=4$, $\sigma_2=2$, and $\sigma_3=2$. }
    \label{fig:fullmodel}
\end{figure}

We present a bifurcation analysis and accompanying basin of attraction plots for Blue success, using the same definition of the basin size from Section~\ref{section:simplemodel}, in Figure~\ref{fig:fullmodel} for both the full model as in \eqref{eq:fullcombatmodel} and the synchronised approximation in \eqref{eq:fullmodelapprox}. We set all cross couplings similarly to the previous section: $\xi_{12} = N_1/ d_T^{(12)}$, $\xi_{13} = N_1/ d_T^{(13)}$, $\xi_{21} = N_2/ d_T^{(21)}$, $\xi_{23} = N_2/ d_T^{(23)}$, $\xi_{31} = N_3/ d_T^{(31)}$, and $\xi_{32} = N_3/ d_T^{(32)}$.

In our scenario, the natural fixed point that corresponds to a Blue victory is $\text{FP}_1 = \left(0,0,K_3(r_3-x_3)/r_3 ,\Delta_1^*,\Delta_2^*\right)$, where Red has been defeated in the competition, Blue has withdrawn due to its success and the centroids have settled to some values $\Delta_1^*$ and $\Delta_2^*$. However, since the $\text{FP}_1$ fixed point is always unstable due to the eigenvalue $\lambda = r_2$, Blue will be unable to achieve victory as Red will persistently surge back from near infinitesimal values. In reality, even in
ecology as in competition, there is
an effective or functional extinction threshold for population resources 
below which it ceases to be viable.
We set here $P_D = 10^{-4}$ as the threshold for a resource below which (i.e. $P_2(t) < P_D $ at some $t$) it is considered defeated and the opposing population successful. In the following simulations, the basins of attraction for a Blue survival are simulated under the constraints of this extinction threshold. 

The results are shown in Figure~\ref{fig:fullmodel} where we show both the linerisation from \eqref{eq:fullmodelapprox} and for the original model in \eqref{eq:fullcombatmodel}. In the latter case, we select high coupling $\sigma_1 = 4$, $\sigma_2 = 2$, and $\sigma_3 = 2$ to demonstrate the validity of our linearisation results when the populations are likely to be synchronised. In these inset plots, the same bifurcation lines derived from the linearisation in \eqref{eq:fullmodelapprox} are superimposed for comparison.
An increase in complexity of observed dynamics is immediately apparent from the panels (a)-(c), where regions of bi-stability and Hopf bifurcations can be observed, notably the exotic looping structures of unstable fixed points in panel (c). In panel (a), a consequence of the imposed death threshold $P_D $ can be observed through the abrupt termination of continuation from the larger valued Hopf point, due to limit cycles beyond that hitting the $P_D$ threshold. As Blue considers victory to be the trivial fixed point, regions in the panels where bi-stability and stable fixed points exist in the panels (a)-(c) will lead to stalemate cases for a proportion of the initial conditions. In panel (b) and directly below in panel (e), we see how the termination of bi-stability around $\phi \approx 0.5$ leads to a sharp transition of the value of the basin size.  In contrast to the previous case in Figure~\ref{fig:simplemodel} where the sign of $\mathcal{K}$ was necessary for the existence of a fixed point, that is we required $|\mu| < 1$, the addition of a third population extends the range of $\mu$ for which we can observe stable fixed points as shown in panel (c). Consequentially, we explore an expanded range of $\mu \in [-\pi,\pi]$. We also note that the $\beta_1$ here is for the dimensional version of the model and that the non-dimensional version can be compared through the transformations $\tilde{\beta_1} = K_1 \beta_1$, $\tilde{\tau} = K_2/K_1 \tau$. A consequence of the refuge provisioning inspired interaction is that Blue is required to adopt a larger value of $\beta_1$ in order to compensate for the refuge provisioning effect because of the shift $(\beta_1 -\beta_1^{(\min)})P_3$, compared with the intrinsic value of $\beta_1$ that would be observed when Green is absent. In contrast to the previous case in Figure~\ref{fig:simplemodel} panels (e) and (f) where an optimal basin size of approximately one can be obtained by making $\beta_1$ large enough, in Figure~\ref{fig:fullmodel} optimal success is also dependent upon superior decision making through an appropriate combination of $\phi$ and $\mu$. Throughout the basin plots in panels (d)-(f), the impact of a third population in the decision making process and the introduction of the death threshold $P_D$, results in additional variation in basin values outside of the boundaries defined by the loss of a non-trivial stable fixed point - being the typical requirement for an optimal basin value. We note however that in panels (e) and (f), the structure of the boundary of basin transition from zero to at least 0.15, as denoted by the dashed black contours, retains a similar shape to the transitions in the earlier figure. 

\section{Parameter significance under the full networked model}\label{sec:numericalstudy}

Due to the complexity of the full model in \eqref{eq:fullcombatmodel} with three populations, we now explore a method to investigate the sensitivity of the basin of attraction of Blue victory to variations within a prescribed subset of the model parameters. We consider this investigation in two parts. The first presents the experimental design, which describes our sampling procedure across the high dimensional parameter space of the model where sampled data is obtained through simulating our model. The second part covers the methods of data analysis using the data obtained from the experimental design. Here we treat the data obtained from the sampling as being separate from the model from which it is simulated. In this analysis, we consider setting the synchronisation feedback function to the power function $p(x) = x^n, \, n \in \mathbb{Z}$ in \eqref{eq:fullcombatmodel} as opposed to the original linear feedback $p(x) = x$. That is we consider the feedback forms $\mathcal{O}_{S_i}^n$ and $\mathcal{O}_{T_i}^n$. We will consider both $n=1$, as originally prescribed in \eqref{eq:fullcombatmodel}, and a second case for $n=2$ which 
corresponds to a sharper modulation as decision making varies between incoherence and synchronisation.

\subsection{Experimental design}
To efficiently explore the vast
parameter space for the size of the
basin of attraction, we use
a hybrid of nearly orthogonal Latin hypercubes (NOLH) \cite{cioppa2007efficient,lin2009construction} and  
Bayesian optimisation \cite{snoek2012practical,brochu2010tutorial}. The aim is to achieve uniform stratification sampling across the basin of attraction for Blue victory. 

The NOLH method \cite{cioppa2007efficient,lin2009construction} allows for the construction of Latin hypercubes, given a fixed sample size, which dramatically improves the space-filling properties of the resultant Latin hypercubes at the expense of inducing small correlations between the columns in the design matrix.  Although NOLH design allows for an efficient approach to exploring parameter space, it does not take the distribution of the basin of attraction values, or `basin value', into consideration. A consequence of this is that the distribution of the basin values will not necessarily be uniform and at worst the distribution might potentially be skewed to the extent that our design `misses' what might be small regions in parameter space, in which Blue's success undergoes significant variation. This is where the Bayesian optimisation procedure
can improve the sampling by intensifying search in parameter regions where prior samples detect relevant
variability. On the other hand, if we had only used Bayesian optimisation, we would still require an initial set of points in our parameter space; randomly sampling is a poor alternative to using a NOLH design in terms of space filling properties.

\begin{algorithm}[tp]
  \caption{Targeted parameter space exploration using Bayesian Optimisation with NOLH initialisation}\label{alg:hybrid}
  \begin{algorithmic}[1]
     \STATE Initialise parameter set using NOLH: $\bm{x} = \{x_i: i=1,\dots,k\} $
     \STATE Initial model evaluation: $\bm{y} = \{y_i = g(x_i): i=1,\dots,k\}$ 
     \STATE Initialise objective function data: $\mathcal{D}_k = \bigcup_{i=1}^k (x_i,f_i(\bm{y}))$
     \STATE Specify Bayesian acquisition function $\alpha$
     \FOR{$n=k+1,2\ldots$}
     \STATE Bayesian acquisition: $\mathbf{x}_n = \text{argmax}_{\mathbf{x}}\alpha(\mathbf{x}|\mathcal{D}_{n-1})$
     \STATE Model evaluation: $y_n = g(\mathbf{x}_n) $
     \STATE Reevaluate objective function and data:
     \FOR{$m=1,2\ldots,n$}
        \STATE  $z_m = f_m(\bm{y})$ 
        \STATE $\mathcal{D}_{m} = \mathcal{D}_{m-1} \bigcup ( x_m,z_m )$
    \ENDFOR
     \STATE Update statistical model
     \ENDFOR
  \end{algorithmic} \label{alg:NOLHbayes}
\end{algorithm}

Thus, in our hybrid approach we first use NOLH to initialise our search. Bayesian optimisation is then used to explore the parameter space, such that we tend towards a uniform stratification sampling across the basin values for Blue victory.  To this end we set the objective function to target a uniform distribution of the basin values which we refer to the as the basin distribution.  A notable difference in this method compared with traditional Bayesian optimisation \cite{shahriari2015taking} is that the objective function $f$, which is a measure of fitness describing the current basin distribution relative to the target basin distribution, is reevaluated at each step. At each new step the value of the basin is obtained and the basin distribution updated for feedback into the objective function $f$. Here we denote the basin value as some $y_i = g(x_i)$ for the parameter set $x_i \in \mathbb{R}^d$ where $d$ is the size of our parameter space and $g$ is the basin valuation function $g:\mathbb{R}^d \rightarrow \mathbb{R}$. This method assumes that the evaluation of $f$ is computationally simple and that the computationally intensive function is $g$. 

In Algorithm~\ref{alg:NOLHbayes}, the objective function for the Bayesian optimisation is given by 
\begin{align} \label{eq:objectivefunction}
    f_m(\bm{y}) =  F(y_m)/ \max_{y\in\bm{y}}\left\{F(y)\right\},
\end{align}
where the function $F(y_m)$ is defined as 
\begin{align}
    F(y_m) = d^2(y_m)/\hat{g}_h(y_m).
\end{align}
In $F(y_m)$, $\hat{g}_h(y_m)$ is the kernel density estimator of the distribution of $\bm{y}$ for some bandwidth $h> 0 $ and $d(y_m)$ is the desired density of $y_m$ where we set $d(y_m) = 1$ to target a uniform distribution. Thus the intent of this procedure is to sample values from our multidimensional parameter space such that we promote uniform stratified sampling of basin values across [0,1] as the objective function \eqref{eq:objectivefunction} promotes sampling from regions in parameter space where there are insufficient samples relative to the target (uniform) sample distribution.

\subsection{Feature selection}

\begin{figure}
    \centering
    \includegraphics[width=\textwidth]{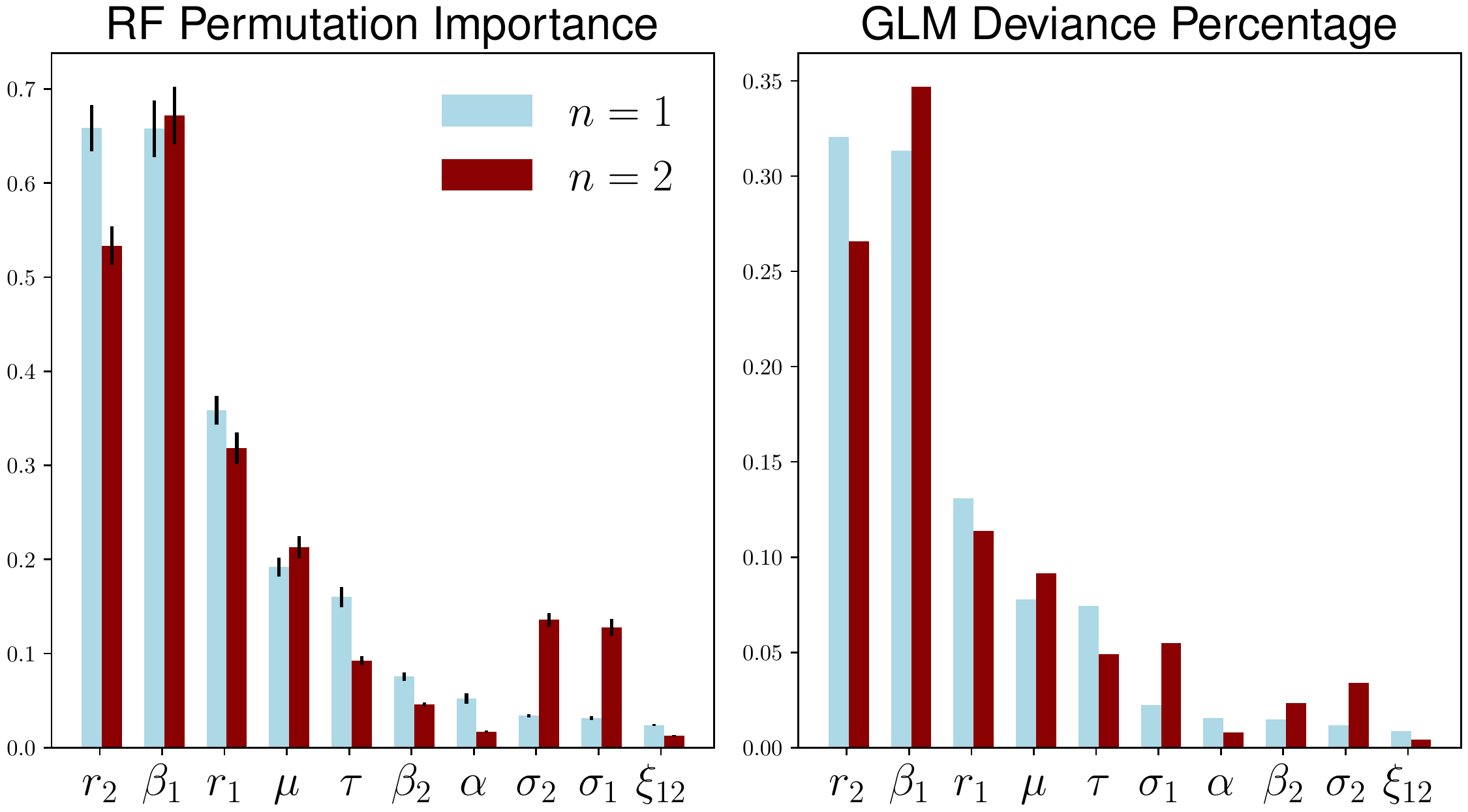}
    \caption{Ranked features through the RF permutation importance (left) and GLM ANOVA deviance percentage (right) for both $n=1$ and $n=2$ versions of the full model in \eqref{eq:fullcombatmodel}. }
    \label{fig:features}
\end{figure}

We now turn to the question of `What are the most important competition processes that lead to a Blue victory?' Using the data set obtained through our experimental design in the previous section, we will identify feature importance through two different methods. We highlight that this analysis will proceed without knowledge of the specific mathematical model from which the data was generated. Our first method consists of a quasi-binomial GLM \cite{mccullagh2019generalized} for the basin size of Blue victory. We note that this approach is simplistic in that it implements parameters as linearly dependent variables without any factoring for their interactions. To thus complement this simple GLM approach, we also consider the application of a Random Forest (RF) regressor \cite{breiman2001random} as a `blackbox' approach to feature selection. Here the RF regressor is implemented through Scikit Learn \cite{pedregosa2011scikit,garreta2013learning} and consists of an averaging over a collection of decision classification trees. For each decision tree a sampled data set is prescribed. The tree then breaks its sampled data set into smaller subsets as it splits nodes. Variable selection for splitting at each node is decided through selecting the variable whose split gives the greatest reduction through a mean square error criterion. We thus typically see variables with the highest importance forming the top levels of the trees ( i.e $\beta_1$, $r_2$, and $r_1$ as in Figure~\ref{fig:features} ), as splitting for these variable results in the largest reduction of variation.

Due to high cardinality in all features, to avoid bias \cite{strobl2007bias} we use permutation importance, defined as the decrease in model score when a feature is randomly shuffled  \cite{breiman2001random}, as our metric for ranking features. The full GLM results are given in the Supplementary material. For the purpose of feature ranking through the GLM, we focus upon the residual deviance from the GLM. Our particular metric is the reduction in deviance attributed to the given parameter from the analysis of deviance table of the GLM. We then take this as a percentage of the total reduction in deviance when comparing the full model with the null model. This percentage provides a comparison with the output from the RF permutation importance. 

In Figure~\ref{fig:features} we plot the results, for $n=1$ and $n=2$, of feature importance for the parameter set obtained using the RF regressor and the deviance percentage from the GLM. To put the following results in context, we recall that Blue suffers a key asymmetry compared to Red in having poorer connectivity
in its network structure. Across the two approaches, the competition parameters $\beta_1$ (Blue's tactical agility in enhancing its capacity to reduce Red), and $r_2$, (Red's recruitment rate) are observed to be the two most significant - these should be regarded as the most important elements to enhance or mitigate to achieve success. Unsurprisingly, for otherwise peer competitors (in size) {\it but with Blue's poorer connectivity}, the capability to reduce the other remains the dominant factor, even in the presence of options for enhancing decision making (through coupling). Evidently Blue could get more out of its decision making capabilities with a more connected structure. However, Figure~\ref{fig:features} now permits understanding of the relative importance of the cognitive dimension. For example, decision making speed, captured in the average frequency difference $\mu$, figures in the middle of the importance sequence. Coupling figures later in the order.

As expected, the relative importance of the internal couplings $\sigma_1$ and $\sigma_2$ notably increase in ranking as we switch between the versions of feedback from $n=1$ to $n=2$. Thus, the more sensitive the  mechanisms for finding and reducing the competitor to the decision making processes, the more value can be obtained from internal coupling. This increasing importance of coupling through increasing $n$ can be alternately achieved by enhancing the network structure; improved connectivity sharpens the transition from low to high values of order parameters. Figure~\ref{fig:features} shows that the more mechanisms available for enhancing interactivity, cohesiveness and agility within the decision making mechanisms, the greater balance against the raw effect of means for reducing, consuming or attriting the competitor.

\section{Discussion and Conclusion}\label{sec:conclusion}

We have developed a hybrid dynamical systems model incorporating distributed decision making within a multi-party adversarial dynamics. Using the Kuramoto model as the basis for the decision making allowed for the role of speed, strength of interaction and patterns of interaction to be explored as counter-balancing factors against the raw factors of reductive strength that otherwise determines the outcomes of competition. With the additional complexity of a third population, using non-trophic mechanisms from mathematical ecology, we were able to see how the presence of such actors create additional demands, in either decision making or reductive power, for each of the adversaries to prevail.

Despite the richness of nonlinear interactions we were able to use a linearisation approximation of the decision dynamics, within an assumption that all actors in the dynamic seek to be, if possible, internally synchronised in their cognitive processes.  This allowed for a fixed point analysis of versions of the general model. Through this, we explored the role of selected decision making processes and competitive agility using bifurcation analysis and numerical simulations. In a simple case, we were able to bound parametric thresholds for competition success; in this sense we described the trade-offs between decision making characteristics and the magnitude of reductive capabilities. Moreover, we characterised the transitions to success in terms of mathematical concepts such as instability of fixed points and basin of attraction size.

For the general model where analytical methods fail, we presented an approach to identify parameters as significant predictors for competition success. Within the simplistic use-case we employed a hierarchically structured group competing against a more randomly, but highly connected, adversary. Indeed we found that reductive power, described through the direct competition parameters, are the leading factors over most of the parameters describing decision making capabilities. However we could explicitly see the impact of increased sensitivity to synchronised decision making on counter-balancing such otherwise raw competitive capabilities. Improved connectivity, correspondingly, allows for the cognitive capabilities to have greater impact on success in competition.  This alone explains why terms such as `decision superiority' figure heavily in  capability in both civilian and military contexts.

One of the main challenges in the construction of such a model lies in the choice of mutual feedback between decision making and competition dynamics. As a new modelling concept, the degree of emphasis on the significance of decision making will vary amongst different modelling practitioners. Though we have used specific expressions for illustrative purposes, the approach that we presented should be viewed as a tool-set for practitioners to test how their assumptions may be realised through various analytical choices for the feedback, and in turn to see the relative importance of the feedback on competition outcome. Thus the present model offers a rich exploratory framework across the range of application areas we outlined from the outset, both in the ecology of species of varying decision making capacities, and in the diverse human endeavours where competition is manifested.


\section*{Acknowledgements}
This research was a collaboration between the Commonwealth of Australia (represented by the Defence Science and Technology Group) and the University of New South Wales Canberra through a Defence Science Partnerships agreement under the auspices of the Modelling Complex Warfighting initiative.

\bibliographystyle{siamplain}
\bibliography{references}
\end{document}


\maketitle

\section{Cluster ansatz for fixed point analysis} \label{supp:clusterthms}

As the system in (4.1)-(4.3) can, in general, only be solved numerically, we explore here an approach to linearise the system under the assumption that each population are internally synchronised.

\subsection{The two-cluster ansatz}

We begin by first defining the centroid difference of the decision cycles between the $P_1$ (Blue) and $P_2$ (Red) populations as 
\begin{align}
\Delta = \overline{\theta}_{1}-\overline{\theta}_{2}
\end{align}
where the centroids are defined as the averages
\begin{align}
\overline{\theta}_{1} = \frac{1}{N_1} \sum_{i\in V_1} \theta_{i},\quad \overline{\theta}_{2} = \frac{1}{N_2} \sum_{j\in V_2} \theta_{j}.
\end{align}

We now repeat the following proposition from the main script:

\begin{proposition1}
Given two internally synchronised populations with phase dynamics described by (1.3) with the centroid difference
\begin{align}
\Delta = \overline{\theta}_{1}-\overline{\theta}_{2},
\end{align}
then a first-order approximation of the phase dynamics in (1.3) is given by 
\begin{align}
    \dot{\Delta} = \mu + S \cos \Delta - C \sin \Delta,
\end{align}
where
\begin{align}
 C & \equiv \frac{d_T^{(12)} \xi_{12} H_1 \cos \phi }{N_1} + \frac{d_T^{(21)} \xi_{21} H_2\cos \psi }{N_2}, \label{eq:Capp} \\
  S & \equiv \frac{d_T^{(12)} \xi_{12} H_1 \sin \phi }{N_1} - \frac{d_T^{(21)} \xi_{21} H_2\sin \psi }{N_2},  \label{eq:Sapp}
\end{align}
with $H_1 = 1-P_2$, $H_2 = 1-P_1$, and $d_T^{(\cdot)} = \sum_{i \in V_1} d_{i}^{(\cdot)}$. 
\end{proposition1}

Since the two populations $V_1$ and $V_2$ are assumed to be synchronised, we consider the ansatz 
\begin{align}
    \theta_i = \overline{\theta}_{1} + b_i,\quad  i \in V_1, \quad 
    \theta_j = \overline{\theta}_{2} + r_j, \quad j \in V_2,
\end{align}
where $b_i$ and $r_j$ are small such that $b_i^2 \approx 0$ and $r_j^2 \approx 0$.
We then linearise the system, keeping only the first term in the Taylor expansion, to obtain
\begin{align}
    \dot{\bm{v}}  + \dot{\bm{\overline{\theta}}} = \mathbf{\Omega} - \mathcal{L} \bm{v} ,\label{eq:linearsys}
\end{align}
where
\begin{align}
    \bm{v} = \left( \begin{array}{c}
         \bm{b}   \\
        \bm{r}  
    \end{array} \right), \quad \dot{\bm{\overline{\theta}}} = \left( \begin{array}{c}
         \overline{\theta}_{1} \bm{1}_{V_1}   \\
        \overline{\theta}_{2} \bm{1}_{V_2} 
    \end{array} \right), \quad \bm{\Omega} = \left( \begin{array}{c}
         \bm{\omega}_1 +\xi_{1} H_1 \sin(\phi - \Delta) \bm{d}^{(12)}    \\
         \bm{\omega}_2 +\xi_{3} H_2 \sin(\psi+ \Delta) \bm{d}^{(21)} 
    \end{array} \right),
\end{align}
and 
\begin{align}
 \mathcal{L} = \left( \begin{array}{cc}
     \sigma_1 H_1 \bm{L}^{(11)} + \xi_{1} H_1\cos(\phi - \Delta) \bm{D}^{(12)} &  -\xi_{1} H_1\cos(\phi - \Delta) \bm{A}^{(12)}\\
     -\xi_{3} H_2\cos(\psi + \Delta)\bm{A}^{(21)} & \sigma_2 H_2 \bm{L}^{(22)} + \xi_{3} H_2\cos(\psi + \Delta) \bm{D}^{(21)}
 \end{array} \right).
\end{align}
In the graph Laplacians within $\mathcal{L}$ and the vector $\mathbf{\Omega}$, the degree of node $i$ within the Blue network connected to Blue nodes and the degree of node $i$ in the Blue network connected to Red nodes are respectively given by
\begin{align}
    d_i^{(11)} = \sum_{j \in V_1} \mathcal{A}_{ij}^{(11)}, \quad d_{ij}^{(12)} = \sum_{j \in V_2} \mathcal{A}_{ij}^{(12)}. 
\end{align}
A similar set of equations $d_i^{(22)}$ and $d_{ij}^{(21)}$ hold for the Red nodes. We then define the degree matrices $\mathbf{D}^{(11)},\mathbf{D}^{(12)}$ and Laplacian $\mathbf{L}^{(11)}$ for Blue through the following equations for their elements: 
\begin{align}
    D_{ij}^{(11)} = d_i^{(11)}\delta_{ij}, \quad D_{ij}^{(12)} = d_i^{(12)}\delta_{ij}, \quad L_{ij}^{(11)} = D_{ij}^{(11)} - \mathcal{A}_{ij}^{(11)}. 
\end{align}
A similar set of equation hold for Red degree matrices and the Laplacian.

Similar to our previous work, we explore the free system basis of the system in \eqref{eq:linearsys}. Here the Blue and Red network Laplacians obey the eigenvalue equations
\begin{align}
    \sum_{j \in V_1} L_{ij}^{(11)} e_j^{(1,r)} &= \lambda_r^{(1)} e_{i}^{(1,r)}, \quad i \in V_1,\\
    \sum_{j \in V_2} L_{ij}^{(22)} e_j^{(2,r)} &= \lambda_r^{(2)} e_{i}^{2,r}, \quad i \in V_2.
\end{align}

As the zeroth eigenvalues $\lambda_0^{(1)}=\lambda_0^{(2)}=0$, the zero mode eigenvectors are given by $\frac{\bm{1}_{V_1}}{N} $ and $\frac{\bm{1}_{V_2}}{M}$. Notably we can state the vector $\dot{\bm{\overline{\theta}}}$ in terms of the zero eigenvectors and we can expand the fluctuations $\dot{\bm{v}}$ through a linear combination of the other eigenvectors - thus the centroids $\overline{\theta}_{1}$ and $\overline{\theta}_{2}$ are the zero-mode projections of the phases $\theta_i, \, i \in V_1$ and $\theta_j, \, j \in V_2$ respectively. Considering the quantity defined from the zero eigenvectors: 

\begin{align}
    \bm{\Pi} = \left( \begin{array}{c}
         \frac{\bm{1}_{V_1}}{N_1}  \\
         -\frac{\bm{1}_{V_2}}{N_2}
    \end{array}\right). \label{eq:projection}
\end{align}

Thus the dynamical system for the centroid difference is given through the following projection:

\begin{align}
    \dot{\Delta} = \bm{\Pi} \cdot \left( \bm{\Omega} - \bm{\mathcal{L}} \bm{v} \right).
\end{align}
Taking the approximation $\dot{\Delta} \approx \bm{\Pi} \cdot \bm{\Omega}$, which is not always guaranteed to hold, yields
\begin{align}
    \dot{\Delta} = \mu + S \cos \Delta - C \sin \Delta,
\end{align}
where
\begin{align}
 \mu & = \overline{\omega}_1 - \overline{\omega}_2, \\
 C & \equiv \frac{d_T^{(12)} \xi_{12} H_1 \cos \phi }{N_1} + \frac{d_T^{(21)} \xi_{21} H_2 \cos \psi }{N_2}, \\
  S & \equiv \frac{d_T^{(12)} \xi_{12} H_1 \sin \phi }{N_1} - \frac{d_T^{(21)} \xi_{21} H_2 \sin \psi }{N_2}, 
\end{align}
and $d_T^{(12)} = \sum_{i \in V_1} d_{i}^{(12)}$.

\begin{proposition2} \label{prop:sup_centroiddynamicsapp}

For some fixed $P_1^*,P_2^* \in [0,1]$, the solution to the centroid dynamics for two internally synchronised populations is given
\begin{align} 
\Delta(t) =  2 \tan^{-1} \left\{  \frac{C-\sqrt{{\cal K} }\tanh \left( \frac{\sqrt{{\cal K}}}{2} (t+\textrm{const} )  \right) }{\mu-S} \right\},
\end{align}
where 
\begin{eqnarray}
{\cal K} \equiv C^2 + S^2 -\mu^2. 
\end{eqnarray}
Furthermore, the existence of a fixed point requires that 
\begin{align}
    C^2 + S^2 \geq \mu^2.
\end{align}
\end{proposition2}

Given some fixed $P_1^*,P_2^* \in [0,1]$, by Proposition~3.1 the dynamics of the centroid are given by
\begin{align}
    \dot{\Delta} = \mu + S \cos \Delta - C \sin \Delta,
\end{align}
where $C$ and $S$ are defined by \eqref{eq:Capp} and ~\eqref{eq:Sapp} in Proposition~3.1 and $P_1$ and $P_2$ are now considered fixed. Using the Weierstrass substitution 
$\Delta = 2 \arctan (\eta)$ where $-\pi < \eta < \pi$ we obtain the new ODE:
\begin{align}
    \dot{\eta} = \frac{1}{2} (\mu-S) \eta^2 -C\eta + \frac{1}{2} (\mu+S)
\end{align}
which we can solve for by using integration by separation and the identity $\tanh ( i \arctan (x)) = ix$ for $-\pi < x < \pi$ to give
\begin{align}
    \eta = \frac{C-\sqrt{{\cal K} }\tanh \left( \frac{\sqrt{{\cal K}}}{2} (t+\textrm{const} )  \right) }{\mu-S},
\end{align}
where 
\begin{eqnarray}
{\cal K} \equiv C^2 + S^2 -\mu^2. \label{eq:Appspecial-K}
\end{eqnarray}

Thus the solution for $\Delta$ is given by 
\begin{align} 
\Delta(t) =  2 \tan^{-1} \left\{  \frac{C-\sqrt{{\cal K} }\tanh \left( \frac{\sqrt{{\cal K}}}{2} (t+\textrm{const} )  \right) }{\mu-S} \right\}
\label{eq:appDeltasol}
\end{align}

It follows that through \eqref{eq:appDeltasol} the sign of $\cal{K}$, where $\cal{K}$ is defined in \eqref{eq:Appspecial-K}, governs the behaviour of the dynamics where if ${\cal K}  < 0$,  $\Delta$ is periodic, when ${\cal K} = 0$ has an unstable fixed point and when ${\cal K}  > 0 $ the exists a stable fixed point.

\subsection{The three-cluster ansatz}

\begin{proposition3} \label{theorem:sup_threecluster_app}
Given three internally synchronised populations with phase dynamics described by (5.2) and the centroid differences
\begin{align}
    \Delta_1 = \overline{\theta}_{1} - \overline{\theta}_{2}, \quad \Delta_2 = \overline{\theta}_{1}- \overline{\theta}_{3}, \quad \Delta_3 = \overline{\theta}_{2} - \overline{\theta}_{3},
\end{align}
then the first-order approximation for synchronised phase dynamics in (5.2) are given by
\begin{align}
    \dot{\Delta}_1 & = \mu - H_1 \left( (\xi_{12} d_T^{(12)}/N_1) \sin\left( \Delta_1 -\phi \right) - (\xi_{13} d_T^{(13)}/N_1) \sin\left( \Delta_2 \right) \right)  \\
    & -H_1  \left( (\xi_{21} d_T^{(21)}/N_2) \sin\left( \Delta_1 +\psi \right) + (\xi_{23}  d_T^{(23)}/N_2) \sin\left( \Delta_2 - \Delta_1\right) \right), \nonumber \\
    \dot{\Delta}_2 & = \nu -H_2 \left(  (\xi_{12} d_T^{(12)}/N_1) \sin\left( \Delta_1 -\phi \right) - \xi_{13} d_T^{(13)}/N_1) \sin\left( \Delta_2 \right) \right) \\
    &  - (\xi_{31} d_T^{(31)}/N_3) \sin\left( \Delta_2 \right) - (\xi_{32} d_T^{(32)}/N_3) \sin\left( \Delta_2 - \Delta_1\right). \nonumber
\end{align}
where $\Delta_3 = \Delta_2 - \Delta_1$, $H_1 = (1-P_2)$, and $H_2 = (1-P_1)$.
\end{proposition3}

Consider the centroids
\begin{align}
    \Delta_1 = \overline{\theta}_{1} - \overline{\theta}_{2}, \quad \Delta_2 = \overline{\theta}_{1} - \overline{\theta}_{3}, \quad \Delta_3 = \overline{\theta}_{2} - \overline{\theta}_{3}
\end{align}

\begin{align}
\overline{\theta}_{1} = \frac{1}{N_1} \sum_{i\in V_1} \theta_{i},\quad \overline{\theta}_{2} = \frac{1}{N_2} \sum_{j\in V_2} \theta_{j}, \quad
\overline{\theta}_{3} = \frac{1}{N_3} \sum_{k\in V_3} \theta_{k}
\end{align}

We explore the fixed point given by the ansatz 
\begin{align}
    \theta_i = \overline{\theta_{1}} + b_i,\quad  i \in V_1, \quad 
    \theta_j = \overline{\theta_{2}} + r_j, \quad j \in V_2, \quad
    \theta_k = \overline{\theta_{3}} + g_j, \quad k \in V_3,
\end{align}
where $b_i$, $r_j$, and $g_k$ are small such that $b_i^2,r_j^2,g_k^2 \approx 0$. The $3 \times 3$ system in \eqref{eq:linearsys} becomes
\begin{align}
    \bm{v} &= \left( \begin{array}{c}
         \bm{b}   \\
        \bm{r}  \\
        \bm{g}
    \end{array} \right), \quad \dot{\bm{\overline{\theta}}} = \left( \begin{array}{c}
         \overline{\theta}_{1} \bm{1}_{V_1}   \\
        \overline{\theta}_{2} \bm{1}_{V_2}  \\
        \overline{\theta}_{3} \bm{1}_{V_3}
    \end{array} \right), \\
    \bm{\Omega} &= \left( \begin{array}{c}
         \bm{\omega}_1+ H_1 \left[ \xi_{12} \sin(\phi - \Delta_1) \bm{d}^{(12)} - \xi_{13} \sin(  \Delta_2) \bm{d}^{(13)} \right]   \\
         \bm{\omega}_2 + H_2  \left[ \xi_{21} \sin(\psi+ \Delta_1) \bm{d}^{(21)} - \xi_{23} \sin(\Delta_3) \bm{d}^{(23)} \right] \\
         \bm{\omega}_3 +\xi_{31} \sin(\Delta_2) \bm{d}^{(31)} + +\xi_{32} \sin(\Delta_3) \bm{d}^{(32)}
    \end{array} \right),
\end{align}
and 

{\small
\begin{align}
 \mathcal{L} = \left( \begin{array}{ccc}
     \sigma_1 H_1 \bm{L}^{(11)} + \mathcal{V}_1  &  -\xi_{12} H_1\cos(\phi - \Delta_1) \bm{A}^{(12)} &  -\xi_{13} H_1\cos( \Delta_1) \bm{A}^{(13)} \\
     -\xi_{21} H_2\cos(\psi + \Delta_1)\bm{A}^{(21)} & \sigma_2 H_2 \bm{L}^{(22)} + \mathcal{V}_2 & -\xi_{23} H_2\cos(\Delta_3)\bm{A}^{(23)} \\
    -\xi_{31} \cos(\Delta_2)\bm{A}^{(31)} &-\xi_{32} \cos(\Delta_3)\bm{A}^{(32)} & \sigma_3 \bm{L}^{(33)} + \mathcal{V}_3
 \end{array} \right)
\end{align}
}%

\begin{align}
    \mathcal{V}_1 &= \xi_{12} H_1\cos(\phi - \Delta_1) \bm{D}^{(12)} + \xi_{13} H_1\cos(\Delta_2) \bm{D}^{(13)}, \\ 
    \mathcal{V}_2 &=  \xi_{21}H_2\cos(\psi + \Delta_1) \bm{D}^{(21)} + \xi_{23} H_2\cos(\Delta_3) \bm{D}^{(23)}, \\ 
    \mathcal{V}_3 &=  \xi_{31} \cos( \Delta_2) \bm{D}^{(31)} + \xi_{32} \cos(\Delta_3) \bm{D}^{(32)}.
\end{align}

Making the approximation that $\dot{\bm{\overline{\theta}}} \approx \bm{\Omega}$ when the fluctuations are considered small,  
the differential equations for the centroids are then given by

\begin{align}
    \dot{\bar{\theta}}_{1} & = \bar{\omega}_1- \frac{\xi_{12} H_1d_T^{(12)}}{N_1} \sin\left( \Delta_1 -\phi \right) - \frac{\xi_{13} H_1d_T^{(13)}}{N_1} \sin\left( \Delta_2 \right)  \\
    \dot{\bar{\theta}}_{2} & = \bar{\omega}_2+ \frac{\xi_{21} H_2d_T^{(21)}}{N_2} \sin\left( \Delta_1 +\psi \right) - \frac{\xi_{23} H_2 d_T^{(23)}}{N_2} \sin\left( \Delta_3 \right)\\
    \dot{\bar{\theta}}_{3} & = \bar{\omega}_3 + \frac{\xi_{31} d_T^{(31)}}{N_3} \sin\left( \Delta_2 \right) + \frac{\xi_{32} d_T^{(32)}}{N_3} \sin\left( \Delta_3 \right)
\end{align}
Since $\Delta_3 = \Delta_2 - \Delta_1$, we can reduce this system to the following two differential equations:

\begin{align}
    \dot{\Delta}_1 & = \mu - H_1 \left( \frac{\xi_{12} d_T^{(12)}}{N_1} \sin\left( \Delta_1 -\phi \right) + \frac{\xi_{13} d_T^{(13)}}{N_1} \sin\left( \Delta_2 \right) \right)  \\
    & - H_2 \left( \frac{\xi_{21} d_T^{(21)}}{N_2} \sin\left( \Delta_1 +\psi \right) - \frac{\xi_{23}  d_T^{(23)}}{N_2} \sin\left( \Delta_2 - \Delta_1\right) \right), \\
    \dot{\Delta}_2 & = \nu - H_1 \left(  \frac{\xi_{12} d_T^{(12)}}{N_1} \sin\left( \Delta_1 -\phi \right) + \frac{\xi_{13} d_T^{(13)}}{N_1} \sin\left( \Delta_2 \right) \right) \\
    &  - \frac{\xi_{31} d_T^{(31)}}{N_3} \sin\left( \Delta_2 \right) - \frac{\xi_{32} d_T^{(32)}}{N_3} \sin\left( \Delta_2 - \Delta_1\right).
\end{align}

\section{Numerical considerations} \label{App:numerical}
\subsection{Simulation procedure}

As is typical in predator-prey dynamics, the outcome of the competition is sensitive to initial conditions; for example, in the
Lanchester setting well-known conservation laws make this dependence explicit. 
While appropriate when excluding decision
making dynamics, this can create
an artificial dependence on the initial
decision state that is a modelling artefact. Note that the steady state of Kuramoto dynamics is indeed independent
of initial conditions. 
Therefore, we initially run the phase dynamics separately with fixed $H_i =1$ for all three populations. This stage of the competition may be regarded as a  reconnaissance period where populations only engage in decision making but have not yet commenced competition, thus preserving their resources up to this point.

In a parameter regime where synchronisation occurs, the phases will settle to a fixed point after an initial transient period. Contrastingly, there are parameter values where phases cannot synchronise, so that the phase dynamics are not necessarily able to settle before we consider simulating the dynamics of the full model after the prescribed transient period. 
The two possibilities are partly contingent on random parameters. In light of this, we compute averages over simulations with $N_{\text{sim}}=100$ where for each instance we randomise the initial phase configurations $\theta_i(t=0)$ for all $i$ but keep the initial population sizes fixed. Thus for each simulation, we fix 
initial conditions for the populations sizes $P_1(t=0), P_2(t=0)$, and $P_3(t=0)$ and average the result of the competition for randomised $\theta_i(0) \in [0,2\pi]$.

Consideration also needs to be made regarding the calculation of the centroids of the phases, $\overline{\theta}_i$. As the decision states measured through $\theta_i$ are equivalent modulo $2\pi$, directly taking the arithmetic mean for the centroid can yield incorrect results if the winding number $q$ is non-zero (from \eqref{eq:winding}). 

\subsection{Evaluation of the centroid of phase oscillators on a circle} \label{app:centroids}

Due to the capacity for individual decision markers to be far ahead of others in their decision loops, the calculation of the centroid $\overline{\theta}$, when phases are considered equivalent modulo $2 \pi$, cannot be made directly through the arithmetic mean when the winding number $q\neq 0$. Here we define  winding number $q\in \mathbb{Z}$ as:
\begin{align}
    q = \frac{1}{2\pi} \sum_{i=1}^N \Delta_{i}, \quad \Delta_{i} =  \left\{ \begin{array}{ll}
    \Delta \theta_i & \text{ for } \quad |\Delta \theta_i| < \pi \\
    \Delta \theta_i - 2\pi & \text{ for } \quad \Delta \theta_i < -\pi \\
    \Delta \theta_i + 2 \pi & \text{ for } \quad  \Delta \theta_i > \pi \end{array} \right. .\label{eq:winding}
\end{align}
where $\Delta \theta_i = \theta_{i+1} - \theta_i$. For a two oscillator system, our solution is to consider an equivalent phase $\theta_2 = \theta_2 + 2\pi k$ for some $k\in \mathbb{Z}$ such that $| \theta_2-\theta_1| < \pi$. Thus under these phase values, we have $q=0$ for which we can calculate the centroid through the arithmetic mean. For a more general system, we consider a similar procedure using a running mean, where at each step we find a $k\in \mathbb{Z}$ such that the winding number between the current running value of the centroid and the new oscillator shifted by $2\pi k$ is $q=0$. That is we calculate:
\begin{align} 
\begin{split}
    \overline{\theta}_1 &= \theta_1 \\
    \overline{\theta}_n &= \overline{\theta}_{n-1} + \frac{1}{n} \left(\theta_n- \overline{\theta}_{n-1} + 2\pi \cdot \argmin_{k_n \in  \mathbb{Z}} \left| \theta_n- \overline{\theta}_{n-1} + 2\pi k_n \right|  \right), \quad n \geq 2.
    \end{split}
\end{align}

Our implementation of this running mean calculation for $\overline{\theta}$ is given in Algorithm~\ref{alg:centroid} below.

\begin{algorithm}[H]
  \caption{Calculation of the centroid of phase oscillators on a circle}\label{alg:centroid}
  \begin{algorithmic}[1]
    \STATE Given $\theta = \left\{\theta_1,\theta_2 ,\dots, \theta_N  \right\}$. 
    \STATE Set $\tilde{\theta} = \mod (\theta, 2\pi)$ and $\theta_C = \tilde{\theta_1}$
     \FOR{$n=2,\ldots,N$}
      \IF{$\|\tilde{\theta} -\theta_C\| > \pi$ }
         \STATE  $\tilde{\theta}_n =  \theta_n + \text{sgn}(\theta_C - \pi) \cdot 2\pi$
      \ELSE
        \STATE  $\tilde{\theta}_n =  \theta_n$
      \ENDIF
      \STATE Update $\theta_C =  \theta_C+ (\tilde{\theta}_n-\theta_C)/n$ 
     \ENDFOR
     \STATE Return $\theta_C$
  \end{algorithmic}
\end{algorithm}

\section{Absent or fixed Green}

When Green ($P_3$) is absent or fixed in its population size and decision making, we can consider the following reduced and nondimensional version of (5.2) as described by the following system: 

\begin{equation}
\begin{aligned} \label{eq:ecomodel}
    \dot{P}_1 & =   r_1 \left(\frac{\alpha P_2 }{1 + \alpha P_2} \right)  P_1 \left(1 - P_1\right ) \mathcal{O}_{S_1}- \beta_2 P_1 P_2 \cdot  \mathcal{O}_{T_2} \frac{\sin(\overline{\theta}_{2}-\overline{\theta}_{1}) +2}{2}- x_1 P_1, \\
    \dot{P}_2 & =  r_2  P_2 \left(1 -P_2 \right ) \mathcal{O}_{S_2} - \frac{\beta_1 P_2}{1+ \tau \beta_1 P_2  } P_1 \cdot \mathcal{O}_{T_1}  \frac{\sin(\overline{\theta}_{1}-\overline{\theta}_{2}) +2}{2}, \\
    \dot{\theta}_i &= \omega_i + H_i \sum_{j\in V}  W_{ij} \sin (\theta_j - \theta_i + \Phi_{ij}) , \quad i \in V. 
    \end{aligned} 
\end{equation}

Although we will use the case where Green is absent, the fixed Green case can be similarly considered through a re-scaling of $\beta_1$. As in the previous model, assuming that the two competing populations are internally synchronised ($\mathcal{O}_1=\mathcal{O}_2=1$), we can apply Proposition~5.1 to reduce the set of ODEs for $\theta$ in \eqref{eq:ecomodel} to the following ODE for the centroid difference:

\begin{align}
 \dot{\Delta}  = \mu + S \cos \Delta - C \sin \Delta, \label{eq:ecomodelapprox}
\end{align}
where $\Delta = \overline{\theta}_{1}-\overline{\theta}_{2}$.
As before, if we consider some fixed $\Delta^* \in \mathbb{R}$ then a set of fixed points can be solved for the reduced system similar to our previous work.
Expressed in terms of Red ($P_2$) , we have five fixed points as follows:

\begin{equation}
\begin{aligned} \label{eq:analyticalFP}
    \text{FP}_1: P_2^* &= 0, \quad \text{FP}_2: P_2^*  = 1, \\
    \text{FP}_3: P_2^* &=  \frac{1}{6 \beta_1 \tau \alpha  r_1 r_2} \left( 2 F_1 + 2^{4/3}\frac{F_2}{F_3}+ 2^{2/3} F_3 \right) , \\
    \text{FP}_4: P_2^* & = \frac{1}{12 \beta_1 \tau \alpha  r_1 r_2} \left( 4 F_1 - 2^{4/3}(1+i\sqrt{3})\frac{F_2}{F_3}- 2^{2/3}(1-i\sqrt{3}) F_3 \right), \\
    \text{FP}_5: P_2^* & =  \frac{1}{12 \beta_1 \tau \alpha  r_1 r_2} \left( 4 F_1 - 2^{4/3}(1-i\sqrt{3})\frac{F_2}{F_3}- 2^{2/3}(1+i\sqrt{3}) F_3 \right),
\end{aligned}
\end{equation}
where the terms $F_1, F_2, F_3 \in \mathbb{C}$, are given by
\begin{align} 
    F_1 &= \alpha  \left(\beta_1\delta \tilde{\beta_2} +(\beta_1 \tau-1) r_1 r_2\right), \\
    F_2 & = \alpha  \left( \alpha \beta_1^2 \delta^2  \tilde{\beta_2} ^2+\beta_1 \delta \xi r_1 r_2+\alpha  r_1^2 r_2 \left(-3 \beta_1^2\delta \tau+\left(1+\beta_1 \tau+\beta_1^2 \tau^2\right) r_2\right)\right), \label{eq:F2} \\
     F_3 & =  \Biggl( \chi  +  \sqrt{ \chi^2-4  \left[F_1^2 +3 \alpha \beta_1 \tau r_1 r_2 \left(\beta_1 \delta (x_1 \alpha +\tilde{\beta_2} )+\alpha  r_1 \left(r_2-\beta_1\right)\right)\right]^3  }
  \Biggr)^{1/3}. \label{eq:F3}
\end{align}
Here we set 
$\xi = (2 \alpha +3) \beta_1  \tilde{\beta_2}  \tau -2 \alpha  \tilde{\beta_2} +3 \alpha  \beta_1  \tau  x_1 $ and 
\begin{align*}
    \chi & = \alpha ^2 \Biggl(2 \beta_1^3 \alpha \delta^3 \tilde{\beta_2} ^3+3 \beta_1^2 \delta^2 \tilde{\beta_2}  \xi r_1 r_2  
     + (\beta_1 \tau-1) \alpha  r_1^3 r_2^2 \left(-9 \beta_1^2\delta \tau+\left(2+5 \beta_1 \tau+2 \beta_1^2 \tau^2\right) r_2\right)   \\
    &    + 3 \beta_1 \delta r_1^2 r_2 \left[-3 \alpha \beta_1^2 \delta \tau   \tilde{\beta_2} + r_2 \left( -\xi+\beta_1  \tau  \left(  \xi + 3\alpha \tilde{\beta_2} + 9 \beta_1\tau x_1  \right)      \right) \right]\Biggr),
\end{align*}
and use the shorthand notation $\tilde{\beta_2} = \beta_2 \cdot \frac{\sin(-\Delta^*) +2}{2} $ and $\delta = \frac{\sin(\Delta^*) +2}{2} $.

The corresponding Blue fixed points $P_1^*$ can be obtained by back substitution into the equation:
\begin{align}
    P_1^* = \frac{r_2}{\beta_1 \delta}(1-P_2^*) (1+\tau \beta_1 P_2^*).
\end{align}

\begin{figure}
    \centering
    \includegraphics[width=\textwidth]{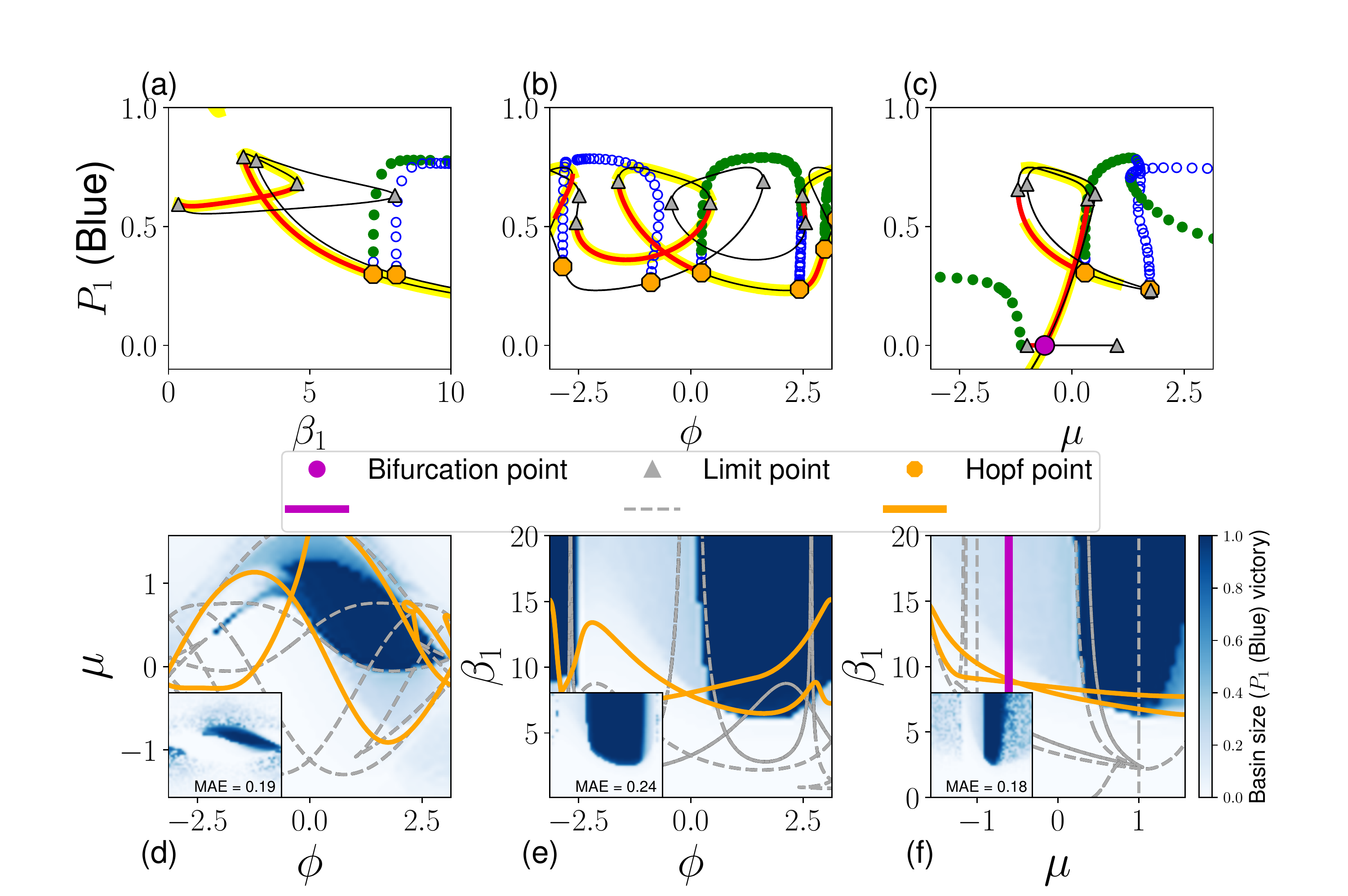}
    \caption{Single parameter (a)-(c) and two parameter (d)-(f) bifurcation diagrams for the extended competitive Kuramoto model with no/fixed Green. In the single parameter plots, red lines denote stable fixed points, black lines denote unstable fixed points, open blue cycles denote unstable limit cycles, and solid green circles denote stable limit cycles. Yellow lines denote the fixed points obtained from solving \eqref{eq:analyticalFP}.  For the two parameter (d)-(f) bifurcation diagrams, we include the basin of attraction size for Blue success. Inset plots show the basin of attraction for the full model for comparison.
    (a): $\phi = 0.2$, $\mu = 0.5$, (b): $\beta_1 = 7.5$, $\mu = 0.25$, (c): $\beta_1 = 7.5$, $\phi = 0.2$, (d): $\beta_1 = 7.5$, (e): $\mu = 0.25$, (f): $\phi = 0.2$}
    \label{fig:eco_model}
\end{figure}

With the fixed point result for $\Delta^*$ from the (3.12), we have a set of equations to be solved for each of the fixed points $\text{FP}_3$,$\text{FP}_4$, and $\text{FP}_5$. Using fixed point iteration, we compare these non-trivial fixed points in Figure~\ref{fig:eco_model} to those obtained through numerical continuation. For the fixed points $\text{FP}_1=(0,0,\Delta^*)$, $\text{FP}_2=(0,1,\Delta^*)$, which correspond to Blue and Red success respectively, we can similarly use our result for the centroid fixed point, where $\Delta^*$ in $\text{FP}_1$ and $\text{FP}_2$ are follow from (3.12).  Given the Jacobian of this system\\

{\footnotesize
\begin{align*}
   \left(
\begin{array}{ccc}
 \frac{2 \alpha  (1-2 P_1) r_1 P_2+(\alpha  P_2+1) \left(\beta_2  P_2 (\sin \Delta-2)-2
   x_1\right)}{2 \alpha  P_2+2} & \frac{1}{2} P_1 \beta_2  (\sin \Delta-2)-\frac{\alpha  (P_1-1)
   P_1 r_1}{(\alpha  P_2+1)^2} & \frac{1}{2} \beta_2 P_1 P_2 \cos \Delta \\
 -\frac{\beta_1  P_2 (\sin \Delta+2)}{2 \beta_1  P_2 \tau +2} & r_2 (1-2 P_2)-\frac{\beta_1  P_1 (\sin
   \Delta+2)}{2 (\beta_1  P_2 \tau +1)^2} & -\frac{\beta_1  P_1 P_2 \cos \Delta}{2 \beta_1  P_2 \tau +2}
   \\
\gamma _2 \sin (\psi + \Delta) & \gamma _1 \sin ( \Delta - \phi ) & -S \sin \Delta - C \cos \Delta \\
\end{array}
\right)
\end{align*}
} %
\\
we can analyse the eigenvalues of the trivial fixed points $E_0$ and $E_1$,
giving 

\begin{align}
\begin{split}
    E_1 &= \quad \lambda_{11} = r_2, \; \lambda_{12} =  -x_1, \; \lambda_{13} =  -\gamma_1\cos(\phi - \Delta^*) -\gamma_2\cos(\psi + \Delta^*). \\
    E_2 &= \quad  \lambda_{21} = -r_2, \; \lambda_{22} =  \beta_2 \frac{(\sin \Delta^* -2)}{2} + \frac{\alpha r_1}{1+\alpha} -x_1 , \; \lambda_{23} = -\gamma_2\cos(\psi + \Delta^*).
    \end{split} \label{eq:ecoeigenvalues}
\end{align}

From \eqref{eq:ecoeigenvalues}, the requirements for stability of the $\text{FP}_1$ and $\text{FP}_2$ fixed points are:
\begin{align*}
  \lambda_{11}: & \quad r_2 < 0; \quad \gamma_1\cos(\phi - \Delta^*) + \gamma_2\cos(\psi + \Delta^*) > 0, \\
  \lambda_{22}: & \quad  \beta_2 > \frac{\frac{\alpha r_1}{1+\alpha} - x_1}{1-\frac{1}{2} \sin \Delta^*}; \quad -\frac{\pi}{2} +\pi k \, \leq \, \psi+\Delta^* \, \leq \, \frac{\pi}{2} + \pi k, \quad k \in \mathbb{Z}.
  \end{align*}

As $r_2>0$ by definition, it follows that the $\text{FP}_1$ fixed point is always unstable and thus Blue will be unable to achieve victory as Red will persistently surge back from near infinitesimal values. In reality, even in
ecology as in conflict there is
an effective or functional extinction threshold for a population 
below which it ceases to be viable.
We denote this here as $P_D = 10^{-4}$ where if a given population falls below this threshold (i.e. $R(t) < P_D$ at some $t$) it is considered defeated and the opposing population successful. In the following simulations, the basins of attraction for a Blue survival are simulated under the constraints of this extinction threshold. 

With the nondimensional parameterisation $\gamma_1 = 1$, $\gamma_2= 1$, $\psi=0$, $\beta_2 =2$, $r_1 =3$, $r_2 =2.5$, $\alpha=20$,  $\tau = 1$, and $x_1=0.25$ we present a set of numerical simulations in Figure~\ref{fig:eco_model}. These results are presented in a similar manner as in previous figures in the main body however the difference in the implementation of the $\beta$ parameter must be noted.
A notable qualitative behaviour of this extended model is the existence of limit cycles arising from the Hopf points (even when $\mathcal{K} >0$ from the decision cycle dynamics) which arises as a consequence of the adaptive form of Blue recruitment through $r_1^*$. An important observation under this parameterisation from the basin of attraction plots in (e) and (f), occurring for both the synchronised approximation of the model \eqref{eq:ecomodelapprox} and the full model \eqref{eq:ecomodel}, is inability of Blue to achieve optimal success through improving $\beta$ if it has an ineffective strategy (i.e $-\pi < \phi < 0 $) or its average collective native frequency is less than Red $\mu < 0$. This is in contrast to the previous model where optimal success can be obtained by making $\beta$ large enough regardless of the other parameters.

\newpage

\section[GLM output]{GLM output}

On Table~\ref{tab:glmtable} where we present the parameter estimate of our quasi-binomial GLM model, the standard error of the estimate, the $t$ value for the Wald test that the regression coefficient for that parameter is 0, and the deviance which measures the contribution to the reduction in deviance of the given parameter as a percentage of the reduction in deviance when comparing the full model to the null model. Here we give values for both versions of the model with $n=1$ and $n=2$. 

\begin{table}[H]
    \centering
    \begin{tabular}{|c|c|c|c|c|c|c|r|r|} \hline
  Coeff.  &  \multicolumn{2}{c|}{Estimate} & \multicolumn{2}{c|}{Std. Error} & \multicolumn{2}{c|}{$t$-value} &   \multicolumn{2}{c|}{Deviance}  \\ \hline \hline
  & n=1 &n=2 & n=1 &n=2 & n=1 &n=2 & n=1 &n=2  \\ \hline
$\mathbf{1}$ & -3.659 &	-3.465&	0.163&	0.176&	-22.483	&-19.693& &
   \\
$\beta_1$        &  0.329	& 0.330&	0.009&	0.009&	37.001&	34.825&		31.34\% &	34.71\%
   \\
$\beta_2$     &   -1.797 &-1.724&	0.151&	0.162&	-11.934&	-10.639	&1.47\%&	2.34\% 
 \\
$\tau$      &   -0.256 &	-0.230&	0.015&	0.016&	-16.704	&-14.338&		7.44\% &	4.92\% 
     \\
$r_1$   & 0.469&	0.481&	0.019&	0.021&	24.960&	22.909&		13.08\% & 	11.36\%
     \\
$r_2$  & -0.766 &	-0.711&	0.021&	0.022&	-37.293&	-32.249&	32.04\%	& 26.57\%
   \\
$\alpha$ & 0.073&0.049&	0.007&	0.008&	9.788&	5.907&		1.56\% &	0.78\%
\\
$\sigma_1$   &  0.109&	0.160&	0.010&	0.011&	11.303&	15.007&		2.24\% &	5.50\%
    \\
$\sigma_2$   &  -0.070&	-0.129&	0.009&	0.010&	-7.408&	-12.927&	1.19\%&	3.41\% 
   \\
$\sigma_3$  & 0.006	&0.008&	0.009&	0.010&	0.650&	0.769& 		0.00\% &	0.00\%
     \\
$\xi_{12}$  &   0.126&	0.076&	0.019&	0.021&	6.511&	3.703&		0.87\%&	0.40\%
     \\
$\xi_{13}$  &  0.006&	-0.011&	0.019&	0.021&	0.331&	-0.513& 		0.01\%&	0.04\%
   \\ 
$\xi_{21}$  &   -0.016&	-0.058&	0.019&	0.020&	-0.801&	-2.859& 0.03\%&	0.11\%
     \\
$\xi_{23}$  &  0.011&	0.050&	0.019&	0.020&	0.562&	2.434& 	0.00\%&	0.12\%
    \\ 
$\xi_{31}$  &  0.005&	-0.027&	0.019&	0.020&	0.279&	-1.323	& 		0.00\%&	0.01\%
    \\
$\xi_{32}$  &   -0.050&	-0.028&	0.019&	0.020&	-2.701&	-1.406&			0.21\%&	0.01\%
   \\
$\phi$     &  0.280&	0.220&	0.047&	0.051&	5.976&	4.345&	0.63\%&	0.40\%
   \\
$\psi$     &   -0.076&	-0.099&	0.048&	0.051&	-1.590&	-1.947&	0.10\%&	0.09\%
      \\
$\mu$ &	0.941&	1.045&	0.048&	0.052&	19.575&	20.184&		7.77\%&	9.16\%\\
$\nu$ &	0.030&	0.090&	0.047&	0.050&	0.628&	1.784& 	0.00\%&	0.06\%\\
 \hline
    \end{tabular}
    \caption{Quasi binomial GLM output for the proportional basin size of Blue victory using a subset of model parameters as dependent variables where deviance represents the reduction in deviance through the addition of the given parameter as a percentage of the overall reduction in deviance when using the full model. Here we obtain deviances from an ANOVA. }
    \label{tab:glmtable}
\end{table}

\newpage

\section[Glossary of parameters,variables, and functions]{Glossary}

On Table~\ref{tab:glossary} we present a glossary of the parameters,variables, and functions present throughout the main manuscript.  
\small
\begin{table}[H]
    \centering
    \begin{tabular}{|c|c|} \hline
  Quantity  &  Description  \\ \hline \hline
    $\alpha$ & Competition resource logistics response rate. \\  
    $\mathbf{A}$ & Adjacency matrix.\\
  $\beta_i$ & Competition reduction rate of population from $i$ onto other. \\
  $\beta_i^{(min)}$ & Minimum competition reduction rate from population $i$ onto others. \\ 
  $\Delta$ & Centroid difference between two populations.\\
  $d_{k}^{(ij)}$ & Degree of node $k$ from adjacency matrix $\mathbf{A}^{(ij)}$.\\
    $d_T^{(ij)}$ & Sum of node degrees $d_{k}^{(ij)}$ from adjacency matrix $\mathbf{A}^{(ij)}$.\\
  $E_i$ & Network edges of population $i$. \\
  $F_{ij}$ & Competition function from population $i$ to $j$.\\
  $\text{FP}$ & A fixed point of the system. \\
  $G_i$ & Graph of population $i$. \\
  $H_i$ & Resource-coupling feedback function for population $i$.\\
  $I_i$ &Initiative function of population $i$. \\
  $K_i$ & Resource capability of population $i$. \\
  $L_{i}$ & Logistic/growth function of population $i$.\\ 
  $\mu$ & Mean intrinsic frequency difference between populations $P_1$ and $P_2$.\\
    $\nu$ & Mean intrinsic frequency difference between populations $P_1$ and $P_3$.\\
  $\omega_j$ & Intrinsic decision making speed of a population $j$.\\
  $\mathcal{O}_{(\cdot)}$ & The Kuramoto order parameter.\\
  $p$ & Synchronisation feedback function.\\    
  $P_i$ & Available resources of population $i$. \\
  $P_D$ & Death threshold for populations. \\
  $\mathbf{\Phi}$ & Frustration matrix. \\
  $\phi$ & $P_1$ decision intent. \\
  $\psi$ & $P_2$ decision intent. \\
  $r_i$ & Recruitment/growth rate of population $i$. \\
  $r^{(max)}_i$ & Maximum recruitment/growth rate of population $i$. \\
  $\sigma_i$ & Intra-network coupling of population $i$. \\
  $\tau$ & Competition search and engagement time. \\
  $\theta_j$ & Decision states' of a population agent $j$. \\
  $\overline{\theta}_i$ & Collective decision state (centroid) of a population $i$. \\
  $V_i$ & Network nodes of population $i$. \\
  $\mathbf{W}$ & The weighted adjacency matrix. \\
  $x_{i}$ & Decay rate of population $i$.\\ 
  $x^{(max)}_{i}$ & Maximum decay rate of population $i$.\\ 
  $x^{(min)}_{i}$ & Minimum decay rate of population $i$.\\ 
  $X_{i}$ & Decay function of population $i$.\\   
  $\xi_{ij}$ & Inter-network coupling from population $i$ to population $j$.\\
  \hline
    \end{tabular}
    \caption{Glossary of the parameters,variables, and functions present throughout the main manuscript.}
    \label{tab:glossary}
\end{table}
